\documentclass[onefignum,onetabnum]{siamart171218}

\usepackage{amsfonts}
\usepackage{graphicx}
\usepackage{epstopdf}
\usepackage{subcaption}
\usepackage[textsize=tiny,textwidth=1.05in]{todonotes}
\usepackage{algorithmic}  
\usepackage{algorithm} 
\usepackage[algo2e]{algorithm2e} 
\usepackage{color}

\ifpdf
  \DeclareGraphicsExtensions{.eps,.pdf,.png,.jpg}
\else
  \DeclareGraphicsExtensions{.eps}
\fi

\newsiamremark{remark}{Remark}
\newsiamremark{hypothesis}{Hypothesis}
\crefname{hypothesis}{Hypothesis}{Hypotheses}
\newsiamthm{claim}{Claim}

\headers{Randomized residual-based error estimators}{K. Smetana, O. Zahm and A. T. Patera}

\title{Randomized residual-based error estimators for parametrized equations\thanks{Submitted to the editors on July 27, 2018.
\funding{This work was funded by ONR Grant N00014-17-1-2077 (ATP).}}}

\author{Kathrin Smetana\thanks{University of Twente, Faculty of Electrical Engineering, Mathematics \& Computer Science, Zilverling,
P.O. Box 217, 7500 AE Enschede, The Netherlands (\email{k.smetana@utwente.nl}).}
\and Olivier Zahm\thanks{Univ. Grenoble Alpes, Inria, CNRS, Grenoble INP, LJK, 38000 Grenoble, France (\email{olivier.zahm@inria.fr}).}
\and Anthony T.~Patera\thanks{Department of Mechanical Engineering, Massachusetts Institute of Technology, Cambridge, MA
02139, United States (\email{patera@mit.edu}).}}

\usepackage{amsopn}

\newcommand{\blue}[1]{#1}

\DeclareMathOperator{\Var}{Var}

\newcommand{\U}{\mathfrak{u}}
\newcommand{\RHS}{\mathfrak{f}}
\newcommand{\tol}{tol}
\newcommand{\qquantile}{q\text{-}\mathrm{quantile}}

\ifpdf
\hypersetup{
  pdftitle={Randomized residual-based error estimators for parametrized equations},
  pdfauthor={K. Smetana, O. Zahm and A. T. Patera}
}
\fi

\begin{document}

\maketitle

\begin{abstract}
We propose a randomized a posteriori error estimator for reduced order approximations of parametrized (partial) differential equations. The error estimator has several important properties: the effectivity is close to unity with prescribed lower and upper bounds at specified high probability; the estimator does not require the calculation of stability (coercivity, or inf-sup) constants; the online cost to evaluate the a posteriori error estimator is commensurate with the cost to find the reduced order approximation; the probabilistic bounds extend to many queries with only modest increase  in cost. To build this estimator, we first estimate the norm of the error with a Monte-Carlo estimator using Gaussian random vectors whose covariance is chosen according to the desired error measure, e.g. user-defined norms or quantity of interest. Then, we introduce a dual problem with random right-hand side the solution of which allows us to rewrite the error estimator in terms of the residual of the original equation. In order to have a fast-to-evaluate estimator, model order reduction methods can be used to approximate the random dual solutions. Here, we propose a greedy algorithm that is guided by a scalar quantity of interest depending on the error estimator. Numerical experiments on a multi-parametric Helmholtz problem demonstrate that this strategy yields rather low-dimensional reduced dual spaces. 
\end{abstract}



\section{Introduction}

Many models for engineering applications, life sciences, environmental issues, or finance depend on parameters which account for variation in the material or geometry but also uncertainty in the data. Often the respective applications require low marginal (i.e. per parameter) computational costs. This is for instance the case in ``many query'' settings where we require the computation of the solution of the corresponding parametrized equation for many different parameter values. 
Examples for model order reduction techniques that aim at computationally feasible approximations of such parametrized models are tensor-based methods \cite{Hac12,Nou17} and the reduced basis (RB) method \cite{HeRoSt16,QuMaNe16,Haa17,RoHuPa08,VePrRoPa03}. In order to ensure say functional safety of a structure, certification of such approximation is of high importance. 
Moreover, bounding the approximation error to get a handle on the uncertainty induced by the approximation is crucial when using it in the context of uncertainty quantification. The subject of this paper is thus certification of approximations to parametrized equations via an a posteriori error estimator for a large number of parameter queries. Our method is also well-suited to real-time contexts. \blue{Employing the a posteriori error estimator say within a greedy algorithm to construct the reduced space requires some (minor) modifications, which we will touch on only very briefly in this paper.}

One of the most commonly used error estimators for inf-sup stable problems is the product of the dual norm of the residual and the inverse of the inf-sup constant. While the former can usually be computed rapidly, accurate estimation of the inf-sup constant is in general rather costly. For instance, the Successive Constraint Method (SCM) \cite{Hetal07,Cetal09,Hetal10} computes a parameter-dependent lower bound of the inf-sup constant by employing the successive solution to appropriate linear optimization problems. This procedure is usually computationally demanding and can lead to pessimistic error bounds \cite{HORU18}.

In this paper we introduce a random a posteriori error estimator which does not require the estimation of stability constants.
The error estimator features several other desirable properties. First, it is both reliable and efficient at given high probability and often has an effectivity close to one. 
Secondly, the effectivity can be bounded from below and above at high probability with constants selected by the user, balancing computational costs and desired sharpness of the estimator.
Moreover, the presented framework yields error estimators with respect to user-defined norms, for instance the $L^2$-norm or the $H^1$-norm; the approach also permits error estimation of linear quantities of interest (QoI).
Finally, depending on the desired effectivity the computation of the error estimator is in general only as costly as the computation of the reduced order approximation or even less expensive, which makes our error estimator strategy attractive from a computational viewpoint.

To derive this error estimator, we consider a Gaussian random vector whose covariance matrix is chosen depending on the respective norm or QoI we wish to estimate. Summing the squares of the inner products of $K$ independent copies of that random vector with the approximation error yields an unbiased Monte Carlo estimator. 
Using concentration inequalities, we control the effectivity of the resulting random error estimator with high probability.
This type of random subspace embedding is typically encountered in compressed sensing \cite{donoho2006compressed}.
\blue{The motivation for using these techniques is to create a high-to-low dimensional map which, in high probability, nearly preserves distances and is thus well-suited for norm estimation.}
By exploiting the error-residual relationship we recognize that these inner products equal the inner products of the residual and the dual solutions of $K$ dual problems with random right-hand sides. Approximating the dual problems via projection-based model order reduction yields an a posteriori error estimator of low marginal computation cost. To construct the dual reduced space we introduce a greedy algorithm 
driven by a scalar QoI that assesses how good the fast-to-evaluate a posteriori error estimator approximates the original Monte Carlo estimator. This goal-oriented strategy outperforms standard dual-residual based greedy algorithms or the Proper Orthogonal Decomposition (POD). \blue{We emphasize that the dual reduced space so obtained does generally \emph{not} contain the primal reduced space as a subspace; the intersection can even be empty. Furthermore, the dimension of the dual reduced space can be smaller than the dimension of the primal reduced space.}

Our a posteriori error estimator is inspired by the probabilistic error estimator for the approximation error in the solution of a system of ordinary differential equations introduced in \cite{CaoPet04} by Cao and Petzold. To estimate the norm of the error, they employ the small statistical sample method from Kenney and Laub \cite{KenLau94}, which estimates the norm of a vector by its inner product with a random vector drawn uniformly at random on the unit sphere.
Rewriting that inner product using the error-residual relationship results in an adjoint (or dual) problem with random final time, whose solution is then invoked to estimate the error \cite{CaoPet04}. This approach is extended to ordinary differential equations via a POD by Homescu et al in \cite{HoPeSe07} and differential algebraic equations in \cite{SeHoPe07}. Also, the effect of perturbations in the initial conditions or parameters on the quality of the approximation of the reduced model is investigated \cite{HoPeSe07,SeHoPe07}. 
In our work we extend these concepts to address the general norms of interest within the PDE context, to explicitly address accurate error estimation for any given parameter value within a finite parameter domain, and to address the limit of many queries.

Randomized methods for error estimation are gaining interest in the reduced order modeling community. 
For instance in \cite{BalNou18}, randomized techniques are used to speed-up the computation of the dual norm of the residual used as an error indicator. By exploiting the fact that the residual manifold is included in a low-dimensional subspace, the authors need appeal to only a few random samples when constructing the random subspace embedding. Instead, our approach targets the true error which, in contrast to the residual, is in general not exactly included in a low-dimensional subspace for the problems we have at hand. Therefore, in our approach, we use different techniques and we determine the number of random sample we need via the cardinality of the parameter set on which we wish to estimate the error.
In \cite{JaNoPr16} a probabilistic a posteriori error bound for linear scalar-valued quantities of interest is proposed, with application in sensitivity analysis. Contrary to the method presented in our work, the right-hand side of the dual problem in \cite{JaNoPr16} is the linear functional associated with the QoI and randomization is done by assuming that the parameter is a random variable on the parameter set.
Another application of randomized techniques, in particular randomized numerical linear algebra \cite{HaMaTr11}, to (localized) model order reduction is considered in \cite{BuhSme18}: a reliable and efficient probabilistic a posteriori error estimator for the difference between a finite-dimensional linear operator and its orthogonal projection onto a reduced space is derived; the main idea is to apply the operator to standard Gaussian random vectors and consider the norm of the result. Also in \cite{zahm2016interpolation}, an interpolation of the operator inverse is built via a Frobenius-norm projection and computed efficiently using randomized methods. An error estimator is obtained by measuring the norm of residual multiplied by the interpolation of the operator inverse, used here as a preconditioner.

We note that also the hierarchical error estimator for the RB method presented in \cite{HORU18} does not require the estimation of any stability constants, such as the inf-sup constant. In \cite{HORU18} the error is estimated by the distance between two reduced approximations of different accuracies and the computational costs depend highly on the dimension of the (primal) reduced space and are always higher than the costs for the computation of the RB approximation. In contrast, in our approach, the costs associated with the dual problems, and hence estimator evaluation, are commensurate with the cost associated with the (primal) RB approximation.
Finally, the reduced-order-model error surrogates (ROMES) method introduced in \cite{DroCar15} and the closely related approaches \cite{MaPaLa16,TrCaDu17,MoStSa18} aim at constructing a statistical model for the approximation error. In \cite{DroCar15} the statistical model is learned via stochastic-process data-fit methods from a small number of computed error indicators.

The remainder of this article is organized as follows. In \cref{sec:NormWithGaussian} we derive a randomized a posteriori error estimator that estimates the error for a finite number of parameter values at given high probability. As this error estimator still depends on the high-dimensional solutions of dual problems, \cref{sec:DualApproximation} is devoted to the reduced order approximation of the dual problems and the analysis of the fast-to-evaluate a posteriori error estimator. In \cref{sec:Numerics} we demonstrate several theoretical aspects of the error estimator numerically and finally draw some conclusions in \cref{sec:Conclusion}.

\newpage

\section{Randomized error estimator for parameter-dependent equation}\label{sec:NormWithGaussian}
\subsection{Parameter-dependent equations and error measurement}\label{sec:motivation}

Consider a real-valued\footnote{Throughout the paper we consider real-valued equations: the extension of our method to the case of complex-valued problems is straightforward using the isomorphy $\mathbb{C}=\mathbb{R}^{2}$.} parameter-dependent equation 
\begin{equation}\label{eq:AUB}
 A(\mu)u(\mu)=f(\mu),
\end{equation}
where the parameter $\mu$ belongs to a parameter set $\mathcal{P}\subset \mathbb{R}^{P}$. 
For every queried parameter $\mu \in \mathcal{P}$, $A(\mu)\in\mathbb{R}^{N\times N}$ is an invertible matrix and $f(\mu)\in\mathbb{R}^{N}$. We assume we are given an approximation $\widetilde u(\mu)$ of the solution $u(\mu)$. In this paper, the goal is to estimate the error
$$
 \| u(\mu)- \widetilde u(\mu)\|_\Sigma .
$$
Here, $\|\cdot\|_\Sigma$ is either a norm defined by means of a symmetric positive-definite (SPD) matrix $\Sigma\in\mathbb{R}^{N\times N}$ via $\|v\|_\Sigma^2 = v^T \Sigma v$ for all $v\in\mathbb{R}^N$, or a semi-norm if $\Sigma$ is only symmetric positive semi-definite. We highlight that the framework presented in this paper encompasses the estimation of the error in various different norms or the error in some QoI as will be discussed in the remainder of this subsection; see \cref{tab:choiceOfSigma} for a brief summary.

By choosing $\Sigma = I_N$, the identity matrix of size $N$, $\|\cdot\|_\Sigma$ becomes the canonical norm $\|\cdot\|_2$ of $\mathbb{R}^N$. If problem \cref{eq:AUB} stems from the discretization of a parameter-dependent linear partial differential equation, there is usually a \emph{natural} norm $\|\cdot\|_X$ associated with a Hilbert space of functions $X\subset H^1(D)$ for some spatial domain $D\subset\mathbb{R}^{d}$, $d\in\{1,2,3\}$. In such a case, there exists a discrete Riesz map $R_X\in\mathbb{R}^{N\times N}$ which is a SPD matrix such that $\sqrt{ (\cdot)^T R_X (\cdot) }= \|\cdot\|_X $.
The choice $\Sigma = R_X$ implies $\|\cdot\|_\Sigma = \|\cdot\|_X$, which means that the error is measured with respect to the natural norm of the problem. We may also consider for instance the error in the $L^{2}$-norm by choosing $\Sigma=R_{L^{2}(D)}$, where the discrete Riesz map $R_{L^{2}(D)}$ is chosen such that $ (\cdot)^T R_{L^{2}(D)} (\cdot) = \|\cdot\|_{L^{2}(D)}^{2}$.

In some cases one is not interested in the solution $u(\mu)$ itself but rather in some QoI defined as a linear function of $u(\mu)$, say 
$$
 s(\mu)=L  u(\mu) ~ \in\mathbb{R}^m,
$$
for some $L\in\mathbb{R}^{m \times N}$.
In this situation one would like to estimate the error $\| s(\mu) - L\,\widetilde u(\mu) \|_W$, where $\|\cdot\|_W$ is a given natural norm on $\mathbb{R}^m$ associated with a SPD matrix $R_W$ so that $\|w\|_W^2 = w^TR_W w$ for all $w\in\mathbb{R}^m$. With the choice $\Sigma = L^T R_W L$ we can write
\begin{align*}
 \| u(\mu) - \widetilde u(\mu) \|_\Sigma^2
 &= (u(\mu)-\widetilde u(\mu))^T \big( L^T R_W L \big) (u(\mu)-\widetilde u(\mu)) =  \| s(\mu) - L\widetilde u(\mu) \|_W^2 ,
\end{align*}
so that measuring the error with respect to the norm $\|\cdot\|_\Sigma$ gives the error associated with the QoI. Notice that if $m<N$ the matrix $\Sigma$ is singular and $\|\cdot\|_\Sigma$ is a semi-norm. Finally, consider the scalar-valued QoI given by $s(\mu)= l^T u(\mu)$ where $l\in\mathbb{R}^N$. This corresponds to previous situation with $m=1$ and $L=l^T$. The choice $\Sigma = l\,l^T $ yields $\| u(\mu) - \widetilde u(\mu) \|_\Sigma^2= |s(\mu) - L\widetilde u(\mu) |$, where $|\cdot|$ denotes the absolute value.

\begin{table}\centering
\begin{tabular}{|l|l|}\hline
 Target error & Choice of $\Sigma$  \\\hline
 $\quad \|u(\mu)-\widetilde u(\mu) \|_2$  & $\quad \Sigma = I_N$ \\
 $\quad \|u(\mu)-\widetilde u(\mu) \|_X$  & $\quad \Sigma = R_X$ \\
 $\quad \|u(\mu)-\widetilde u(\mu)\|_{L^{2}(D)}$ & $\quad\Sigma=R_{L^{2}(D)}$\\
 $\quad \|s(\mu)-L\widetilde u(\mu) \|_W$ & $\quad \Sigma = L^T R_W L$ \\
 $\quad |s(\mu)-l^T \widetilde u(\mu) |$  & $\quad \Sigma = l\,l^T$ \\\hline
\end{tabular}
\caption{Possible choices for $\Sigma$ depending on the target error.}
\label{tab:choiceOfSigma}
\vspace{-15pt}
\end{table}

\subsection{Estimating norms using Gaussian maps}\label{subsec:estimateNorm}

\blue{
In this section we show how the (semi-)norm $\|\cdot\|_\Sigma$ can be approximated by $\|\Phi\cdot\|_2$ for some random matrix $\Phi \in\mathbb{R}^{K\times N}$ with $K\ll N$.

Let $Z\sim\mathcal{N}(0,\Sigma)$ be a zero mean Gaussian random vector in $\mathbb{R}^N$ whose covariance matrix is chosen to be the matrix $\Sigma\in\mathbb{R}^{N\times N}$ which defines the (semi-)norm $\|\cdot\|_\Sigma$, cf \cref{tab:choiceOfSigma}}.
Given a vector $v\in\mathbb{R}^N$, for example $v=u(\mu)-\widetilde u(\mu)$ for some (fixed) parameter $\mu\in\mathcal{P}$, we can write
$$
 \|v\|_\Sigma^2 
 = v^T \Sigma v 
 = v^T \mathbb{E}( ZZ^T )v 
 = \mathbb{E}( (Z^T v)^2 ),
$$
where $\mathbb{E}(\cdot)$ denotes the \blue{expected value}.
This means that $(Z^T v)^2$ is an unbiased estimator of $\|v\|_\Sigma^2$. Let $Z_1,\hdots,Z_K$ be $K$ independent copies of $Z$ and define the random matrix $\Phi\in\mathbb{R}^{K\times N}$ whose $i$-th row is $(1/\sqrt{K})Z_i^T$. The matrix $\Phi$ is sometimes called a \textit{Gaussian map}. Denoting by $\|\cdot\|_2$ the canonical norm of $\mathbb{R}^K$, we can write
\begin{equation}\label{eq:defPhi}
 \| \Phi v \|_2^2 = \frac{1}{K}\sum_{i=1}^K  ( Z_i^T v )^2 \qquad \text{for any } v\in \mathbb{R}^N.
\end{equation}
In other words, $\| \Phi v \|_2^2$ is a $K$-sample Monte-Carlo estimator of $\mathbb{E}( (Z^T v )^2 ) = \|v\|_\Sigma^2$. By the independence of the $Z_i$'s, we have $\Var( \| \Phi v \|_2^2  ) = \frac{1}{K} \Var( Z^T v  )$ so that $\| \Phi v \|_2^2$ is a lower variance estimator of $\|v\|_\Sigma^2$ compared to $(Z ^T v )^2$. However, the variance is not always the most relevant criteria to assess the performance of an estimator. In the context of this paper, we rather want to quantify the probability that $\| \Phi v \|_2^2$ deviates from $\|v\|_\Sigma^2$. This can be done by noting that, provided $\|v\|_\Sigma\neq0$, the random variables $(Z_i^T v ) / \|v\|_\Sigma$ for $i=1,\hdots,K$ are independent standard normal random variables so that we have
$$
 \| \Phi v \|_2^2 = \frac{\|v\|_\Sigma^2}{K} \sum_{i=1}^K \Big( \frac{Z_i^T v}{\|v\|_\Sigma}\Big)^2 \sim \frac{\|v\|_\Sigma^2}{K} Q ,
$$
where $Q \sim \chi^2(K)$ follows a chi-squared distribution with $K$ degrees of freedom. Denoting by $\mathbb{P}\{A\}$ the probability of an event $A$ and by $\overline{A}$ the complementary event of $A$, the previous relation yields
\begin{align*}
 \mathbb{P}\Big\{ w^{-1}\|v\|_\Sigma \leq \| \Phi v \|_2 \leq w \|v\|_\Sigma \Big\}
 &=1-\mathbb{P}\big\{  \overline{ Kw^{-2} \leq  Q \leq Kw^2 }  \big\} ,
\end{align*}
for any $w\geq1$. Then for any given (fixed) vector $v\in \mathbb{R}^N$, the probability that a realization of $\| \Phi v \|_2$ lies between $w^{-1}\|v\|_\Sigma$ and $w\|v\|_\Sigma$ is independent of $v$ but also independent of the dimension $N$. The following proposition gives an upper bound for $\mathbb{P}\big\{  \overline{ Kw^{-2} \leq  Q \leq Kw^2 }  \big\} $ in terms of $w$ and $K$. The proof, given in \cref{proof:Chi2Tail}, relies on the fact that we have closed form expressions for the law of $Q\sim\chi^2(K)$.

\begin{proposition}\label{prop:Chi2Tail}
 Let $Q\sim\chi^2(K)$ be a chi-squared random variable with $K\geq3$ degrees of freedom. For any $w>\sqrt{e}$ we have 
 \begin{equation*}\label{eq:Chi2Tail}
  \mathbb{P}\big\{  \overline{ Kw^{-2} \leq  Q \leq Kw^2 }  \big\} 
  \leq \Big( \frac{\sqrt{e}}{w} \Big)^{K} .
 \end{equation*}
\end{proposition}

\begin{figure}[t]
\begin{center}
  \begin{subfigure}[c]{0.47\textwidth}\centering
  \includegraphics[width=\textwidth]{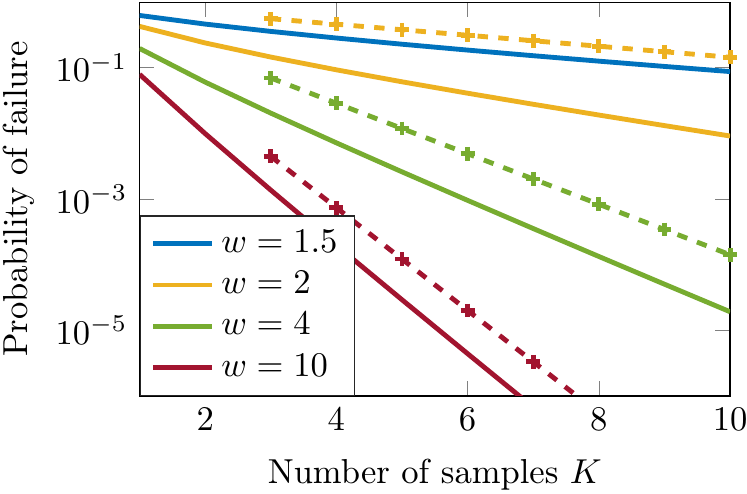}
  \end{subfigure}
  \begin{subfigure}[c]{0.45\textwidth}\centering
   \footnotesize
   \begin{tabular}{|c|cc|} \hline
     $w$ & \multicolumn{2}{|c|}{$K=3$} \\ \hline
     1.1 & $8.2\times 10^{-1}$ & $-$  \\
     2 & $1.4\times 10^{-1}$ & $5.6\times 10^{-1}$  \\
     5 & $1.0\times 10^{-2}$ & $3.5\times 10^{-2}$  \\
     10 & $1.3\times 10^{-3}$ & $4.4\times 10^{-3}$ \\
     50 & $1.1\times 10^{-5}$ & $3.5\times 10^{-5}$ \\\hline \hline
     
     $w$ & \multicolumn{2}{|c|}{$K=10$} \\ \hline
     1.1 & $6.7\times 10^{-1}$ & $-$ \\
     2 & $9.1\times 10^{-3}$ & $1.4\times 10^{-1}$ \\
     5 & $2.2\times 10^{-6}$ & $1.5\times 10^{-5}$ \\
     10 & $2.4\times 10^{-9}$ & $1.4\times 10^{-8}$ \\
     50 & $2.6\times 10^{-16}$ & $1.5\times 10^{-15}$ \\\hline
   \end{tabular}
  \end{subfigure}
\end{center}
\caption{Exact value of $\mathbb{P}\big\{  \overline{ Kw^{-2} \leq  Q \leq Kw^2 }  \big\} $ (solid curves on the graph, left column on the table) and its upper bound $( \sqrt{e}/w)^{K}$ given by \cref{prop:Chi2Tail} (dashed curves on the graph, right column on the table) for different values of $K$ and $w$.}
\label{fig:Chi2}
\vspace{-15pt}
\end{figure}

\cref{prop:Chi2Tail} shows that the probability $\mathbb{P}\big\{  \overline{ Kw^{-2} \leq  Q \leq Kw^2 }  \big\}$ decays at least exponentially with respect to $K$, provided $w\geq\sqrt{e}$ and $K\geq3$. 
Then for any $v\in \mathbb{R}^N$, the relation
\begin{equation}\label{eq:controlError}
 w^{-1}\|v\|_\Sigma \leq \| \Phi v \|_2 \leq w \|v\|_\Sigma ,
\end{equation}
holds with a probability greater than $1-(\sqrt{e}/w)^K$. As expected, a large value of $w$ is beneficial to ensure the probability of failure $(\sqrt{e}/w)^K$ to be small. For instance with $w=4$ and $K=6$, relation \cref{eq:controlError} holds with a probability larger than $0.995$. However, we observe in \cref{fig:Chi2} that this theoretical result is rather pessimistic since it overestimates the true probability by one order of magnitude for small values of $w$. Also, we conjecture on \cref{fig:Chi2} that there is an exponential decay even when $w\leq \sqrt{e}$ (see the blue curve with $w=1.5$), which is not predicted by \cref{prop:Chi2Tail}.

In many situations we want to estimate the norm of several vectors rather than just one vector solely. This is for instance the case if one has to estimate the norm of the error $v= u(\mu) - \widetilde u(\mu)$ for many different parameter values $\mu \in\mathcal{P}$. In that case, one would like to quantify the probability that relation \cref{eq:controlError} holds simultaneously for any vector in a set $\mathcal{M}\subset \mathbb{R}^N$. Assuming $\mathcal{M}$ is finite, a union bound argument
--- for a detailed proof see \cref{proof:estimate many vectors} --- yields the following result:

\begin{corollary}\label{prop:estimate many vectors}
Given a finite collection of vectors $\mathcal{M}=\{v_1,v_2,\hdots, v_{\# \mathcal{M}}\}\subset \blue{\mathbb{R}^N}$ and a failure probability $0<\delta<1$. Then, for any $w > \sqrt{e}$ and 
\begin{equation}\label{eq:ConditionK_FiniteSet}
  K\geq \min \left\lbrace \frac{\log( \# \mathcal{M}) + \log(\delta^{-1})}{\log(w/\sqrt{e})} , \enspace 3 \right \rbrace
\end{equation}
we have
\begin{equation}\label{eq:estimate many vectors}
 \mathbb{P}\Big\{ w^{-1}\|v\|_\Sigma \leq \| \Phi v \|_2  \leq w  \|v\|_\Sigma ~,~\forall v\in\mathcal{M} \Big\} \geq 1 - \delta. 
\end{equation}
\end{corollary}

\cref{tab:concentration_set} gives numerical values of $K$ that satisfy \cref{eq:ConditionK_FiniteSet} depending on $\delta$, $w$ and $\#\mathcal{M}$. For example with $\delta=10^{-4}$ and $w=10$, estimating simultaneously the norm of $10^{9}$ vectors requires only $K=17$ samples. Again, we emphasize that this result is independent on the dimension $N$ of the vectors to be estimated.

\begin{table}[t]
  \centering
  \footnotesize
  
  \begin{tabular}{| l | c c c | c | l | c c c |} 
      \cline{2-4}\cline{7-9}
      \multicolumn{1}{c|}{$\delta=10^{-2}$}& $w=2$ & $w=4$ & $w=10$ &
      \multicolumn{1}{c}{$\quad$}& 
      \multicolumn{1}{c|}{$\delta=10^{-4}$}& $w=2$ & $w=4$ & $w=10$ 
      \\  
      \cline{1-4}\cline{6-9}
      $\#\mathcal{M}=10^0$   & 24 & 6  &  3  &     & $\#\mathcal{M}=10^0$ & 48  & 11 & 6  \\ 
      $\#\mathcal{M}=10^3$   & 60 & 13 &  7  &     & $\#\mathcal{M}=10^3$ & 84  & 19 & 9  \\ 
      $\#\mathcal{M}=10^6$   & 96 & 21 & 11  &     & $\#\mathcal{M}=10^6$ & 120 & 26 & 13  \\ 
      $\#\mathcal{M}=10^9$   & 132& 29 & 15  &     & $\#\mathcal{M}=10^9$ & 155 & 34 & 17  \\ 
      \cline{1-4}\cline{6-9}
  \end{tabular} 
  
 \caption{Minimal value of $K$ for which Condition \cref{eq:ConditionK_FiniteSet} is satisfied.
}
\label{tab:concentration_set}
\vspace{-15pt}
\end{table}

\begin{remark}[Comparison with the Johnson-Lindenstrauss lemma \cite{dasgupta1999elementary,johnson1984extensions}] 
 The Johnson-Lindenstrauss (JL) lemma states that for any $0<\varepsilon<1$ and any finite set $\mathcal{M}\subset\mathbb{R}^N$, the condition $K\geq 8\varepsilon^{-2}\log(\#\mathcal{M})$ ensures the existence of a linear map $\Phi: \mathbb{R}^{N} \rightarrow \mathbb{R}^{K}$ such that
 $
 (1-\varepsilon) \|v-u\|_2^2 \leq \| \Phi v -\Phi u\|_2^2 \leq (1+\varepsilon)\| v-u \|_2^2,
 $
 holds for all $u,v\in\mathcal{M}$. 
 Replacing $\mathcal{M}$ by $\mathcal{M}\cup\{0\}$ and letting $u=0$, one has that $K\geq 8\varepsilon^{-2}\log(\#\mathcal{M}+1)$ is sufficient to ensure the existence of a $\Phi\in\mathbb{R}^{K\times N}$ such that
 \begin{equation}\label{eq:tmp29647}
  \sqrt{1-\varepsilon} \|v\|_2 \leq \| \Phi v \|_2 \leq \sqrt{1+\varepsilon} \| v \|_2, \quad \text{for all } v\in\mathcal{M}.
 \end{equation}
The above relation differs from \cref{eq:controlError} in the sense that the deviation of $ \| \Phi v \|_2$ from $\|v\|_2$ is controlled in an \emph{additive manner} via a parameter $\varepsilon$ instead of a \emph{multiplicative way} via $w$. We highlight also the different dependencies of $K$ on $\varepsilon$ and $w$. In contrast to the requirement in the JL lemma Condition \cref{eq:ConditionK_FiniteSet} permits reduction in the number of required copies $K$ of the random vectors by considering an increased $w$. Note that the computational complexity of the a posteriori error estimator we propose in this paper crucially depends on $K$, see \cref{sec:online complexity}. Since the goal in this paper is to estimate the error we do in general not have to insist on a very accurate estimation of $\|v\|_2$. Instead, in many situtations it might be preferable to accept a higher effectivity $w$ of the a posteriori error estimator in favour of a faster computational time. We emphasize that the user has the choice here. 
 
Notice also that with the choice $w=1/\sqrt{1-\varepsilon}$, Equation \cref{eq:tmp29647} implies \cref{eq:controlError}. Then, the JL lemma ensures that \cref{eq:controlError} holds true if $K\geq 8(1-w^{-2})^{-2}\log(\#\mathcal{M}+1)$. Even if we have the same logarithmic dependence on $\#\mathcal{M}$, this is much larger than what we obtained in \cref{eq:ConditionK_FiniteSet}, already for moderate but especially for large values of $w$. For example with $w=4$, $\#\mathcal{M}=10^3$ and $\delta=10^{-2}$, JL lemma requires $K\geq 63$ whereas Condition \cref{eq:ConditionK_FiniteSet} requires only $K\geq 13$.
Finally, we highlight that a similar result to \cref{eq:controlError} has been obtained in \cite{KenLau94} for random vectors that are uniformly and randomly selected from the sphere. The multiplicative type of estimates in \cite{KenLau94} motivated us to derive similar results for Gaussian vectors. 

 
\end{remark}

\begin{remark}[Drawing Gaussian vectors]
In actual practice we can draw efficiently from $Z\sim\mathcal{N}(0,\Sigma)$ using a factorization of the covariance matrix of the form of $\Sigma=U^{T}U$, \textit{e.g.} a (sparse) Cholesky decomposition.
It is then sufficient to draw a standard Gaussian vector $\widehat{Z}$ and to compute the matrix-vector product $Z=U^{T}\widehat{Z}$.
As pointed-out in \cite[Remark 2.9]{BalNou18}, one can take advantage of \blue{a} potential block structure of $\Sigma$ to build a (non-square) factorization $U$ with a negligible computational cost.
\end{remark}

\subsection{Randomized a posteriori error estimator}\label{subsec:estimate many vectors}

We apply the methodology described in the previous subsection to derive a residual-based randomized a posteriori \blue{error} estimator for the error $\| u(\mu)-\widetilde u(\mu)\|_\Sigma$. Let $\Phi = K^{-1/2}[Z_1,\hdots Z_K]^T$ be a random matrix in $\blue{\mathbb{R}^{K\times N}}$ where $Z_1,\hdots Z_K$ are independent copies of $Z\sim\mathcal{N}(0,\Sigma)$, and consider the error estimator $\Delta(\mu)=\| \Phi \big( u(\mu)-\widetilde u(\mu) \big) \|_2$, or equivalently
\begin{equation}\label{defeq:error est}
 \Delta(\mu)  
 = \left( \frac{1}{K} \sum_{k=1}^{K} \Big(Z_{i}^{T} \big(u(\mu)-\widetilde u(\mu) \big) \Big)^{2} \right)^{1/2}.
\end{equation}
If the parameter set $\mathcal{P}$ is finite, \cref{prop:estimate many vectors} with $\mathcal{M}=\{ u(\mu)-\widetilde u(\mu) ; \mu\in\mathcal{P} \}$ permits control of the quality of the estimate $\Delta(\mu)$ uniformly over $\mu\in\mathcal{P}$. But in actual practice the parameter set is often of infinite cardinality. Using more sophisticated techniques than just a simple union bound argument should provide results also when $\mathcal{P}$ has infinite cardinality. In this paper, we are however only interested in the case of a finite set of parameter values, as restated in the following corollary.

\begin{corollary}\label{coro:truth est S}
 Let $0 < \delta < 1$ and $w > \sqrt{e}$. Given a finite set of parameter values $\mathcal{S}\subset\mathcal{P}$, the condition
 \begin{equation}\label{coro:truth est Condition S}
   K\geq \min \left\lbrace \frac{\log( \# \mathcal{S}) + \log(\delta^{-1})}{\log(w/\sqrt{e})} , \enspace 3 \right \rbrace,
 \end{equation}
 is sufficient to ensure
 $$
  \mathbb{P}\Big\{ w^{-1} \Delta(\mu) \leq \|u(\mu)-\widetilde u(\mu) \|_\Sigma \leq w \Delta(\mu) ~,~\forall \mu \in \mathcal{S} \Big\} \geq 1 - \delta. 
 $$
\end{corollary}

It is important to note that Condition \cref{coro:truth est Condition S} depends only on the cardinality of $\mathcal{S}$. This means that $K$ can be determined only knowing the \emph{number of parameters} for which we need to estimate the error.
However, computing $\Delta(\mu)$ requires the solution $u(\mu)$ of problem \cref{eq:AUB}, which is infeasible in practice. By introducing the residual
\begin{equation}\label{eq:PrimalResidual}
 r(\mu)= f(\mu) - A(\mu)\widetilde u(\mu),
\end{equation}
associated with Problem \cref{eq:AUB} and, similar to \cite{CaoPet04,HoPeSe07}, exploiting the \emph{error residual relationship} we may albeit rewrite the terms $Z_{i}^{T} (u(\mu)- \widetilde u(\mu))$, $1\leq i \leq K$ as follows:
\begin{equation}\label{eq:err res relationship}
Z_{i}^{T} (u(\mu) - \widetilde u(\mu)) = Z_{i}^{T} A(\mu)^{-1} r(\mu) = (A(\mu)^{-T}Z_{i})^{T}r(\mu).
\end{equation}
\blue{The terms $Z_{i}^{T} (u(\mu)- \widetilde u(\mu))$ thus equal the inner products of the (primal) residual and the solutions $Y_{i}(\mu) \in \mathbb{R}^N$} of the \emph{random dual problems} 
\begin{equation}\label{eq:randomDualProblem}
 A(\mu)^T Y_i(\mu) = Z_i, \quad 1\leq i \leq K.
\end{equation}

Because of the random right hand side in \cref{eq:randomDualProblem}, the solutions $Y_1(\mu),\hdots,Y_K(\mu)$ are random vectors. 
Thanks to the above the error estimator $\Delta(\mu)$ \cref{defeq:error est} can be rewritten as
\begin{equation}
 \Delta(\mu) 
 = \left( \frac{1}{K}\sum_{i=1}^K \big( Y_i(\mu)^T r(\mu) \big)^2 \right)^{1/2} . \label{eq:DualTrick}
\end{equation}
This shows that $\Delta(\mu)$ can be computed by applying $K$ linear forms to the residual $r(\mu)$.
In that sense, $\Delta(\mu)$ can be considered as an a posteriori error estimator.
Notice that computing the solutions to \cref{eq:randomDualProblem} is in general as expensive as solving the primal problem \cref{eq:AUB}. In the next section we show how to approximate the dual solutions $Y_1(\mu),\hdots,Y_K(\mu)$ in order to obtain a fast-to-evaluate a posteriori error estimator.

\begin{remark}[Scalar-valued QoI]
 When estimating the error in scalar-valued QoIs of the form of $s(\mu)=l^T u(\mu)$, the covariance matrix is $\Sigma=l\,l^T$, see  \cref{sec:motivation}. In that case the random vector $Z\sim\mathcal{N}(0,\Sigma)$ follows the same distribution as $X\,l$ where $X\sim\mathcal{N}(0,1)$ is a standard normal random variable (scalar). The random dual problem \cref{eq:randomDualProblem} then becomes
 $
  A(\mu)^TY_i(\mu) = X_i\, l 
 $
 and the solution is $Y_i(\mu) = X_i \, q(\mu)$ where $q(\mu)$ is the solution of the deterministic dual problem $A(\mu)^T q(\mu) = l$. Dual problems of this form are commonly encountered for estimating linear quantities of interest, see \cite{pierce2000adjoint} for a general presentation and \cite{Haa17,RoHuPa08,zahm2017projection} for the application in reduced order modeling.
 
\end{remark}

\blue{
\begin{remark}[Considerations when employing $\Delta(\mu)$ to enrich the reduced space]
Say that we use the a posteriori error estimator $\Delta(\mu)$ to select a new parameter and use the associated solution to enrich the reduced space. Then, we wish to use $\Delta(\mu)$ again for the enriched reduced space. However, now the problem occurs that the error between $u(\mu)$ and the reduced solution that uses the newly selected snapshot depends on the error estimator and thus $Z_{1},\hdots,Z_{K}$; we lose independence. One solution would be to redraw the samples in each iteration, which is in general however computationally infeasible. Alternatively, as suggested in \cite{BalNou18}, we can adapt the number of samples $K$ in order to take into account (using union bound arguments) all possible outcomes of the greedy algorithm; for further details we refer to \cite[Section 5.1]{BalNou18}.
\end{remark}}

\section{A fast-to-evaluate randomized a posteriori error estimator}\label{sec:DualApproximation}

In order to obtain a fast-to-evaluate a posteriori error estimator whose computational complexity is independent of $N$, we employ projection-based model order reduction (MOR) techniques to compute approximations of the solutions $Y_1(\mu),\hdots,Y_K(\mu)$ of the dual problems \cref{eq:randomDualProblem}. To that end, let us assume that we are given a \blue{fixed} realization of the $K$ random vectors $Z_1,\hdots,Z_K$ and that we have a reduced space $\widetilde{\mathcal{Y}}\subset \mathbb{R}^N$ at our disposal.
\blue{Different ways to construct $\widetilde{\mathcal{Y}}$ will be discussed in \cref{sec:greedy} and compared numerically in \cref{sec:Numerics}. In any case, $\widetilde{\mathcal{Y}}$ will be built from dual solutions $Y_{i}(\mu)$, $i=1,\hdots,K$ of \eqref{eq:randomDualProblem} for random right-hand sides
$Z_1,\hdots,Z_K$, the latter being fixed before constructing the dual reduced space. $\widetilde{\mathcal{Y}}$ should thus be considered as a random subspace}. 
Then, we define $\widetilde Y_i(\mu)$ as the Galerkin projection of $Y_i(\mu)$ on $\widetilde{\mathcal{Y}}$, meaning
\begin{equation}\label{eq:monolithic}
 \widetilde Y_i(\mu) \in \widetilde{\mathcal{Y}} \,: \quad
 \langle  A(\mu)^T\widetilde Y_i(\mu) ,v \rangle = \langle  Z_i ,v \rangle
 \,, \quad \forall v\in\widetilde{\mathcal{Y}}.
\end{equation}
Here, $\langle v,w\rangle:=v^{T}w$ for all $v,w \in \mathbb{R}^{N}$. 
We emphasize that we employ the same reduced space $\widetilde{\mathcal{Y}}$ for the approximation of the $K$  dual solutions $Y_1(\mu),\hdots,Y_K(\mu)$. Needless to say that a \emph{segregated strategy}, \blue{where we construct and use $K$ different dual reduced spaces $\widetilde{\mathcal{Y}}_i$ for the $K$ different right-hand sides $Z_{i}$ and associated dual reduced solutions $Y_{i}(\mu)$, $i=1,\hdots,K$}, can also be considered. The advantage of a segregated strategy is that one can easily parallelize the computations, if needed. However, in this paper we focus exclusively on the \emph{monolithic approach} \cref{eq:monolithic}\blue{, employing one single dual reduced space}.

By replacing $Y_i(\mu)$ in \cref{eq:DualTrick} by the fast-to-evaluate approximation $\widetilde Y_i(\mu)$, we define a fast-to-evaluate a posteriori error estimator as
\begin{equation}\label{eq: def a post est online}
 \widetilde \Delta(\mu) := \left( \frac{1}{K} \sum_{i=1}^{K} (\widetilde Y_{i}(\mu)^{T} r(\mu))^{2} \right)^{1/2}. 
\end{equation}
We highlight that, in constrast to for instance the ``standard'' a posteriori error estimator being defined as the product of the reciprocal of a stability constant and the dual norm of the primal residual, $\widetilde \Delta(\mu)$ does not contain any constants that require estimation. Moreover, unlike hierarchical error estimators \cite{BMNP04,HORU18} the quality of the approximation used for the error estimator does \emph{not} depend on the quality of the primal approximation\blue{; the dual reduced space does not in general contain the primal reduced space as a subspace and can even be of smaller dimension than the latter.} For a more elaborate comparison we refer to \cref{sec:online complexity}.

Additionally, we shall show in \cref{sec:online complexity} that evaluating $\mu \mapsto \widetilde \Delta(\mu)$ requires only the solution of one linear system of size $n_{\widetilde{\mathcal{Y}}}:=\dim(\widetilde{\mathcal{Y}})$, instead of $K$ linear systems of size $n_{\widetilde{\mathcal{Y}}}$ as suggested by \cref{eq:monolithic}. 
However, before discussing the computational complexity of $\widetilde \Delta(\mu)$, we show in \cref{sec:performance of online-efficient estimator} that under certain conditions $\widetilde \Delta(\mu)$ is both a reliable and efficient error estimator at high probability. Based on this analysis we propose in \cref{sec:greedy} different greedy algorithms for constructing the reduced space $\widetilde{\mathcal{Y}}$.

\subsection{Analysis of the fast-to-evaluate a posteriori error estimator}\label{sec:performance of online-efficient estimator}

First, we relate the relative error in the a posteriori error estimator to the error in the dual residual:

\begin{proposition}\label{prop:dualError}
 Assume $\Sigma$ is invertible. The fast-to-evaluate error estimator $\widetilde\Delta(\mu)$ defined by \cref{eq: def a post est online} satisfies
 \begin{align}
  \frac{| \Delta(\mu) - \widetilde\Delta(\mu) |}{\|u(\mu)-\widetilde u(\mu)\|_\Sigma}
  &\leq \max_{1\leq i\leq K} \| A^T(\mu) \widetilde Y_i(\mu) - Z_i \|_{\Sigma^{-1}}  \qquad \text{ for all }  \mu\in\mathcal{P}.\label{eq:dualError}
 \end{align}
Here, $\|\cdot\|_{\Sigma^{-1}}$ denotes the norm on $\mathbb{R}^{N}$ such that $\|v\|_{\Sigma^{-1}}^2 = v^T\Sigma^{-1}v$ for all $v\in\mathbb{R}^N$.
 
\end{proposition}

The proof is given in \cref{proof:dualError}. Notice that \cref{prop:dualError} requires $\Sigma$ to be invertible, which excludes the cases where one wants to estimate the error in a vector-valued QoI, see \cref{sec:motivation}. \cref{prop:dualError} allows us to control the error via $\widetilde\Delta(\mu)$, where the effectivity $w$ is enlarged in an additive manner, as stated in the following corollary.

\begin{corollary}\label{cor:dualErrorAdditive}
 Suppose we are given a finite set of parameter values $\mathcal{S}\subset\mathcal{P}$ for which we want to estimate the error $\|u(\mu)-\widetilde u(\mu) \|_\Sigma $.
 Let $0 < \delta < 1$, $w > \sqrt{e}$ and assume
 \begin{equation}\label{eq:assumption dualErrorAdditive}
   K\geq \min \left\lbrace \frac{\log( \# \mathcal{S}) + \log(\delta^{-1})}{\log(w/\sqrt{e})} , \enspace 3 \right \rbrace.
 \end{equation}
 Furthermore, assume that $\Sigma$ is invertible and that we have $\varepsilon\leq w^{-1}$, where
 \begin{equation}\label{eq:eTilde2}
  \varepsilon = \sup_{\mu\in\mathcal{P}} \left\{ \max_{1\leq i\leq K} \| A^T(\mu) \widetilde Y_i(\mu) - Z_i \|_{\Sigma^{-1}} \right\} .
 \end{equation}
 Then, we have
 \begin{equation}\label{eq:dualErrorAdditive}
  \mathbb{P}\Big\{ 
  (w + \varepsilon)^{-1} \widetilde\Delta(\mu) 
  \leq \|u(\mu)-\widetilde u(\mu) \|_\Sigma 
  \leq \frac{w}{1 - w \, \varepsilon } \,\widetilde\Delta(\mu) ,
  \quad \forall\mu\in\mathcal{S}
  \Big\} \geq 1-\delta .
 \end{equation}
\end{corollary}

The proof is given in \cref{proof:dualErrorAdditive}.
\cref{cor:dualErrorAdditive} gives a sufficient condition to control the quality of the estimator $\widetilde\Delta(\mu)$ over a finite set of parameter values $\mathcal{S}\subset\mathcal{P}$ with high probability. It requires $\varepsilon\leq w^{-1}$, which is equivalent to $\| A^T(\mu) \widetilde Y_i(\mu) - Z_i \|_{\Sigma^{-1}} \leq w^{-1}$ for all $\mu\in\mathcal{P}$ and all $1\leq i \leq K$. 
\blue{To satisfy this condition, one has to design an algorithm which builds $\widetilde Y_i(\mu)$ in a way that $A^T(\mu) \widetilde Y_i(\mu)$ is close to $Z_i$ uniformly over $\mathcal{S}$ and independently on the value taken by $Z_i$. Obtaining $\varepsilon\leq w^{-1}$ can however be challenging (from a computational perspective).} To explain this, let us note that $\|Z_i\|_{\Sigma^{-1}}$ is, with high probability\footnote{To show this, note that $\|Z_i\|_{\Sigma^{-1}}^2\sim \chi^2(N)$ so that, by \cref{prop:Chi2Tail}, relation $w'^{-1}\sqrt{N}\leq \|Z_i\|_{\Sigma^{-1}}\leq w'\sqrt{N}$ holds with probability $1-(\sqrt{e}/w')^N$ for any $w'\geq\sqrt{e}$.}, of the order of $\sqrt{N}$.
Therefore $\varepsilon\leq w^{-1}$ means that the \textit{relative} dual residual norm ought to be of the order of
$$
 \frac{\| A^T(\mu) \widetilde Y_i(\mu) - Z_i \|_{\Sigma^{-1}}}{ \| Z_i \|_{\Sigma^{-1}} }
 \simeq
 \frac{1}{w\sqrt{N}}.
$$
When $N\gg1$, the condition $\varepsilon\leq w^{-1}$ means that we need a very accurate approximation of the dual variables. For instance with $N=10^6$, the dual residual norm $\| A^T(\mu) \widetilde Y_i(\mu) - Z_i \|_{\Sigma^{-1}} / \| Z_i \|_{\Sigma^{-1}}$ has to be less that $10^{-3}$ for all $\mu\in\mathcal{P}$ and all $1\leq i\leq K$, which can be too demanding in actual practice.

Next, we give an alternative way of controlling the quality of $\widetilde\Delta(\mu)$.
Contrarily to \cref{cor:dualErrorAdditive}, which provides a \textit{additive} type of control, the following proposition gives a control in an \textit{multiplicative} manner. The proof in given in \cref{proof:dualErrorMultiplicative} in the appendix.

\begin{proposition}\label{prop:dualErrorMultiplicative}
 Suppose we are given a finite set of parameter values $\mathcal{S}\subset\mathcal{P}$ over which we want to estimate the error $\|u(\mu)-\widetilde u(\mu) \|_\Sigma $.
 Let $0 < \delta < 1$, $w > \sqrt{e}$ and assume
 \begin{equation}\label{eq:assumption dualErrorMultiplicative}
   K\geq \min \left\lbrace \frac{\log( \# \mathcal{S}) + \log(\delta^{-1})}{\log(w/\sqrt{e})} , \enspace 3 \right \rbrace.
 \end{equation}
 Then the fast-to-evaluate estimator $\widetilde\Delta(\mu)$ satisfies
 \begin{equation}\label{eq:dualErrorMultiplicative}
  \mathbb{P}\Big\{ 
  (\alpha w)^{-1} \widetilde\Delta(\mu) 
  \leq \|u(\mu)-\widetilde u(\mu) \|_\Sigma 
  \leq (\alpha w) \,\widetilde\Delta(\mu),
  \quad \mu \in \mathcal{S},
  \Big\} \geq 1-\delta ,
 \end{equation}
 where
 \begin{equation}\label{eq:alpha}
  \alpha := \, \max_{\mu \in \mathcal{P}} \left( \max \left\{\frac{\Delta(\mu)}{\widetilde \Delta(\mu)} \,,\, \frac{\widetilde \Delta(\mu)}{\Delta(\mu)} \right\} \right) \geq 1.
 \end{equation}
\end{proposition}

\cref{prop:dualErrorMultiplicative} shows that, with high probability, the error estimator $\widetilde\Delta(\mu)$ departs from the true error $\|u(\mu)-\widetilde u(\mu) \|_\Sigma$ at most by a multiplicative factor $(\alpha w)^{-1}$ or $(\alpha w)$. 
Notice that $\alpha$ is a measure of the distance from $\mu\mapsto\widetilde\Delta(\mu)$ to $\mu\mapsto\Delta(\mu)$: if it is close to $1$ then $\widetilde\Delta(\mu)$ is close to $\Delta(\mu)$ uniformly over the parameter set $\mathcal{P}$. Unlike \cref{cor:dualErrorAdditive}, \cref{prop:dualErrorMultiplicative} does not require $\Sigma$ to be invertible and, even more importantly, \blue{it does not put any restrictions on $\alpha$.} However, the computation of $\alpha$ can be expensive since it requires the exact error estimator $\Delta(\mu)$ over the whole parameter set $\mathcal{P}$. Therefore, we propose to use $\alpha$ as a stopping criterion when constructing the dual reduced space to ensure that $(\alpha w)^{-1} \widetilde\Delta(\mu) 
  \leq \|u(\mu)-\widetilde u(\mu) \|_\Sigma 
  \leq (\alpha w) \,\widetilde\Delta(\mu)$ holds true for a rich training set $\subset \mathcal{P}$ as we will detail in \cref{sec:greedy}.

\subsection{Greedy constructions of the dual reduced space $\widetilde{\mathcal{Y}}$}\label{sec:greedy}

\subsubsection{Vector point of view of the dual problems}

A popular technique to build a reduced space is to take the span of snapshots of the solution. In order to handle the $K$ distinct dual problems, the index ``$i$'' in \cref{eq:randomDualProblem} in considered as an additional parameter. Thus, we define the augmented parameter set $\mathcal{P}_K = \{ 1,\hdots,K \} \times \mathcal{P}$ and seek a $n_{\widetilde{\mathcal{Y}}}$-dimensional reduced space of the form of
\begin{equation}\label{eq:RBdualSpaceVector}
 \widetilde{\mathcal{Y}} = \text{span}\{ Y_{i_1}(\mu_1),\hdots,Y_{i_{n_{\widetilde{\mathcal{Y}}}}}(\mu_{n_{\widetilde{\mathcal{Y}}}}) \} ,
\end{equation}
where the $n_{\widetilde{\mathcal{Y}}}$ elements $(i_1,\mu_1),\hdots,(i_{n_{\widetilde{\mathcal{Y}}}},\mu_{n_{\widetilde{\mathcal{Y}}}})$ are to be chosen in $\mathcal{P}_K$. The RB methodology (see for instance \cite{HeRoSt16,QuMaNe16,Haa17,RoHuPa08} for an introduction) consists in selecting $(i_1,\mu_1),\hdots,(i_{n_{\widetilde{\mathcal{Y}}}},\mu_{n_{\widetilde{\mathcal{Y}}}})$ in a greedy fashion \cite{VePrRoPa03}. In detail, assuming \blue{that} the $j$ first parameters are given, the $(j+1)$-th parameter is defined as
\begin{equation}\label{eq:weakGreedy_training}
 (i_{j+1},\mu_{j+1}) \in \underset{ (i,\mu)\in\mathcal{P}_K^\text{train} }{\text{argmax}} \| A(\mu)^T \widetilde Y_i(\mu) - Z_i \|_* ,
\end{equation}
where $\widetilde Y_i(\mu)$ is the approximation of $Y_i(\mu)$ given by \cref{eq:monolithic} with $\widetilde{\mathcal{Y}}$ defined as in \cref{eq:RBdualSpaceVector}. 
Here, $\mathcal{P}_K^\text{train} \subset \mathcal{P}_K$ is a sufficiently rich training set with finite cardinality and $\|\cdot\|_*$ denotes an arbitrary norm of $\mathbb{R}^N$.
According to \cref{cor:dualErrorAdditive}, it is natural to chose $\|\cdot\|_* = \|\cdot\|_{\Sigma^{-1}}$, provided $\Sigma$ is invertible. After having computed the snapshot $Y_{i_{j+1}}(\mu_{j+1})$, the reduced space $\widetilde{\mathcal{Y}}$ is updated using \cref{eq:RBdualSpaceVector} with $j\leftarrow j+1$. By selecting the parameter $(i_{j+1},\mu_{j+1})$ according to \cref{eq:weakGreedy_training}, the idea is to construct a reduced space that minimizes the dual residual norm $\| A(\mu)^T \widetilde Y_i(\mu) - Z_i \|_*$ uniformly over the training set $(i,\mu)\in\mathcal{P}_K^\text{train} $.

It remains to define a criterion to stop the greedy iterations.
Given a user-defined tolerance $\tol\geq0$, the use of the stopping criterion
\begin{equation}\label{eq:stoppingCriterion MAX STD}
 \max_{(i,\mu)\in\mathcal{P}_K^\text{train}} \| A(\mu)^T \widetilde Y_i(\mu) - Z_i \|_* \leq \tol ,
\end{equation}
ensures that, at the end of the iteration procedure, the residual norm of the dual problem is below $\tol$ everywhere on the training set $\mathcal{P}_K^\text{train}$.
One can relax that criterion by replacing the max in \cref{eq:stoppingCriterion MAX STD} by the quantile of order $ q \in[0,1]$:
\begin{equation}\label{eq:stoppingCriterion QUANTILE STD}
 \qquantile \Big\{ \| A(\mu)^T \widetilde Y_i(\mu) - Z_i \|_* \,:\,  (i,\mu)\in\mathcal{P}_K^\text{train} \Big\} \leq \tol .
\end{equation}
\blue{Here $\qquantile\{A\}$ denotes the $\lceil q\#A \rceil$-th largest entries of a (ordered and finite) set $A$.}
With this stopping criterion, the iterations stop when at least a fraction of $q$ points in $\mathcal{P}_K^\text{train}$ have a dual residual norm below $\tol$. 
Notice that the $\qquantile$ and the $\max$ coincides when $q=1$ so that \cref{eq:stoppingCriterion QUANTILE STD} generalizes \cref{eq:stoppingCriterion MAX STD}.
The resulting greedy algorithm is summarized in \cref{algo:greedy}.

\begin{algorithm}[t]
 \caption{Greedy construction of the dual reduced space $\widetilde{\mathcal{Y}}$}\label{algo:greedy}
 \KwData{Operator $\mu\mapsto A(\mu)$, samples $\{Z_1,\hdots,Z_K\}$, training set $\mathcal{P}_K^\text{train}$, tolerance $\tol$, quantile order $q$}
 
 Initialize $\widetilde{\mathcal{Y}}=\{0\}$ and $j=0$
 
 \While{$\qquantile_{(i,\mu)\in\mathcal{P}_K^\text{train}}\{ \| A(\mu)^T \widetilde Y_i(\mu) - Z_i \|_* \} > \tol$}{
  Compute $\widetilde Y_i(\mu) \in \widetilde{\mathcal{Y}}$ via \cref{eq:monolithic} 
  
  Find $(i_{j+1},\mu_{j+1})$ that maximizes $(i,\mu)\mapsto \| A(\mu)^T \widetilde Y_i(\mu) - Z_i  \|_*$ over $\mathcal{P}_K^\text{train}$
  
  Compute the snapshot $Y_{i_{j+1}}(\mu_{j+1}) = A(\mu_{j+1})^{-T}Z_{i_{j+1}}$
  
  Update the dual reduced space $\widetilde{\mathcal{Y}} \leftarrow \widetilde{\mathcal{Y}} + \text{span}\{ Y_{i_{j+1}}(\mu_{j+1}) \}$
  
  Update $j\leftarrow j+1$
 }
 \KwResult{Dual reduced space $\widetilde{\mathcal{Y}}$.}
\end{algorithm}

\subsubsection{Matrix point of view of the dual problems}

{\blue{Next,}} we propose another greedy algorithm which relies on a matrix interpretation of the $K$ dual problems \cref{eq:randomDualProblem}. Let us denote by 
$$
 \mathbf{Y}(\mu) =[ Y_1(\mu),\hdots, Y_K(\mu)] \in\mathbb{R}^{N\times K} ,
$$
the matrix containing the dual solutions. Instead of constructing the reduced space $\widetilde{\mathcal{Y}}$ as the span of vectors $Y_i(\mu)$, like in Equation \cref{eq:RBdualSpaceVector}, we now consider reduced spaces of the form of
\begin{equation}\label{eq:RBdualSpaceMatrix}
 \widetilde{\mathcal{Y}} = \text{span}\{ \mathbf{Y}(\mu_1) \lambda_1 , \hdots, \mathbf{Y}(\mu_{n_{\widetilde{\mathcal{Y}}}}) \lambda_{n_{\widetilde{\mathcal{Y}}}} \},
\end{equation}
where $\lambda_1,\hdots,\lambda_{n_{\widetilde{\mathcal{Y}}}}$ are $n_{\widetilde{\mathcal{Y}}}$ vectors in $\mathbb{R}^K$ and where $\mu_1,\hdots,\mu_{n_{\widetilde{\mathcal{Y}}}} \in\mathcal{P}$.
Notice that if the vectors $\lambda_i$ are canonical vectors of $\mathbb{R}^K$ we have that $\widetilde{\mathcal{Y}}$ can be written as in \cref{eq:RBdualSpaceVector}.
In that sense, the approximation format \cref{eq:RBdualSpaceMatrix} is richer than \cref{eq:RBdualSpaceVector} and we can expect better performance. We also note that the greedy algorithm we propose here shares some similarities with the POD-greedy algorithm introduced in \cite{HaaOhl08}.

We now propose a second greedy algorithm inspired by \cref{prop:dualErrorMultiplicative}. Let $\mathcal{P}^\text{train} \subset \mathcal{P}$ be again a finite training set and suppose that at step $r$ in the greedy algorithm we have a reduced space $\widetilde{\mathcal{Y}}$ as in \cref{eq:RBdualSpaceMatrix} at our disposal. The first step is to define the next evaluation point $\mu_{j+1}$ as
\begin{equation}\label{eq:strongGreedy}
 \mu_{j+1} \in \underset{\mu \in \mathcal{P}^\text{train}}{\text{argmax}} \left( \max \left\{\frac{\Delta(\mu)}{\widetilde \Delta(\mu)} \,,\, \frac{\widetilde \Delta(\mu)}{\Delta(\mu)} \right\} \right) ,
\end{equation}
where we recall that $\Delta(\mu)=(\frac{1}{K}\sum_{i=1}^{K}[Z_{i}^{T}(u(\mu)-\widetilde u(\mu))]^{2})^{1/2}$. 
Finding $\mu_{j+1}$ according to \cref{eq:strongGreedy} requires to compute the solution $u(\mu)$ over the training set $\mu\in\mathcal{P}^\text{train}$. \blue{As this is in general not computationally efficient, we suggest replacing $u(\mu)$ by a reference solution $u_{\mathrm{ref}}(\mu)$ such that $\|u_{\mathrm{ref}}(\mu) - u(\mu)\|_{\Sigma} \ll \| \widetilde u(\mu) - u(\mu)\|_{\Sigma}$. We can choose as $u_{\mathrm{ref}}(\mu)$ for instance a hierarchical approximation of $u(\mu)$, where we use a larger primal reduced space to determine $u_{\mathrm{ref}}(\mu)$. Note that we only suggest using such a reference solution for the construction of the dual reduced space and not afterwards when certifying the reduced approximation in the online stage.}
Then, we introduce the reference error estimator 
\begin{equation*}\label{eq:ref_est}
\Delta_{\mathrm{ref}}(\mu):=\left(\frac{1}{K}\sum_{i=1}^{K}\big( Z_{i}^{T}(u_{\mathrm{ref}}(\mu)-\widetilde u(\mu))\big)^{2}\right)^{1/2}
\end{equation*}
and seek $\mu_{j+1}$ as
\begin{equation}\label{eq:strongGreedy_bis}
 \mu_{r+1} \in \underset{\mu \in \mathcal{P}^\text{train}}{\text{argmax}} \left( \max \left\{\frac{\Delta_{\mathrm{ref}}(\mu)}{\widetilde \Delta(\mu)} \,,\, \frac{\widetilde \Delta(\mu)}{\Delta_{\mathrm{ref}}(\mu)} \right\} \right).
\end{equation}
Once the parameter $\mu_{j+1}$ is found either with \cref{eq:strongGreedy} of with \cref{eq:strongGreedy_bis}, we compute the dual solutions $Y_1(\mu_{j+1}),\hdots,Y_K(\mu_{j+1})$ and assemble $\mathbf{Y}(\mu_{j+1})$. Here we need to solve $K$ linear equations with the same operator $A(\mu_{j+1})^T$ but with $K$ different right-hand sides, see Equation \cref{eq:randomDualProblem}. This can be done efficiently say be using a Cholesky or LU decomposition and reusing the factorization for the $K$ problems.

The second step is to determine the vector $\lambda_{j+1}$. In order to maximize the improvement of the reduced space, we propose to define $\lambda_{j+1}$ as follows:
\begin{equation}\label{eq:strongGreedyLambda}
 \lambda_{j+1} \in \underset{\lambda \in\mathbb{R}^K}{\text{argmax }} \frac{ \| \mathbf{Y}(\mu_{j+1}) \lambda - \widetilde{\mathbf{Y}}(\mu_{j+1}) \lambda   \|_2 }{\|\lambda\|_2 } ,
\end{equation}
where $\widetilde{\mathbf{Y}}(\mu_{j+1}) = [ \widetilde Y_1(\mu_{j+1}),\hdots, \widetilde Y_K(\mu_{j+1})]$. The rational behind \cref{eq:strongGreedyLambda} is to align $\lambda_{j+1}$ with the direction where the matrix $\widetilde{\mathbf{Y}}(\mu_{j+1})$ differs the most from $\mathbf{Y}(\mu_{j+1})$. One can easily show that $\lambda_{j+1}$ defined by \cref{eq:strongGreedyLambda} is the first eigenvector of the $K$-by-$K$ matrix 
\begin{equation}\label{eq:matrixM}
 M(\mu_{j+1}) = \big(  \mathbf{Y}(\mu_{j+1}) - \widetilde{\mathbf{Y}}(\mu_{j+1})  \big)^T\big( \mathbf{Y}(\mu_{j+1}) - \widetilde{\mathbf{Y}}(\mu_{j+1}) \big) .
\end{equation}
Once $\lambda_{j+1}$ is computed, we set $j\leftarrow j+1$ and we update the reduced space $\widetilde{\mathcal{Y}}$ using \cref{eq:RBdualSpaceMatrix}. We terminate the algorithm based on the following stopping criteria
$$
 \qquantile \left\{ \max \left\{\frac{\Delta(\mu)}{\widetilde \Delta(\mu)} \,,\, \frac{\widetilde \Delta(\mu)}{\Delta(\mu)} \right\}
 \,:\, \mu\in\mathcal{P}^\text{train} \right\}
 \leq \tol.
$$
The resulting greedy algorithm is summarized in \cref{algo:qoi_greedy}.

\begin{algorithm}[t]
 \caption{Greedy construction of $\widetilde{\mathcal{Y}}$ with goal oriented greedy selection}\label{algo:qoi_greedy}
 \KwData{Operator $\mu\mapsto A(\mu)$, samples $\{Z_1,\hdots,Z_K\}$, training set $\mathcal{P}^\text{train}$, tolerance $\tol$, quantile order $q$, approximation $\mu\mapsto \widetilde u(\mu)$, reference solution $\mu\mapsto u_{\mathrm{ref}}(\mu)$}
 ~\\[-0.7cm]
 
 Compute $\Delta_{\mathrm{ref}}(\mu)$ for all $\mu\in\mathcal{P}^\text{train}$
 
 Initialize $\widetilde{\mathcal{Y}}=\{0\}$ and $j=0$

 \While{$\qquantile_{\mu\in\mathcal{P}^\text{train}}\big\{ \max \left\{\frac{\Delta_{\mathrm{ref}}(\mu)}{\widetilde \Delta(\mu)} \,,\, \frac{\widetilde \Delta(\mu)}{\Delta_{\mathrm{ref}}(\mu)} \right\} \big\} > \tol$}{
 
  Define $\widetilde Y_i(\mu) \in \widetilde{\mathcal{Y}}$ by \cref{eq:monolithic} and $\widetilde \Delta(\mu)$ by \cref{eq: def a post est online}
  
  Find $\mu_{j+1}$ that maximizes $\mu\mapsto \max \left\{\frac{\Delta_{\mathrm{ref}}(\mu)}{\widetilde \Delta(\mu)} \,,\, \frac{\widetilde \Delta(\mu)}{\Delta_{\mathrm{ref}}(\mu)} \right\}$ over $\mathcal{P}^\text{train}$
  
  Compute the solutions $Y_{i}(\mu_{j+1}) = A(\mu_{j+1})^{-T}Z_{i}$ for all $1\leq i \leq K$ 
  
  Compute the matrix $M(\mu_{j+1})$ by \cref{eq:matrixM} and its leading eigenvector $\lambda_{j+1}$
  
  Update the dual reduced space $\widetilde{\mathcal{Y}} \leftarrow \widetilde{\mathcal{Y}} + \text{span}\{ \mathbf{Y}(\mu_{j+1}) \lambda_{j+1} \}$
  
  Update $j\leftarrow j+1$
 }
 \KwResult{Dual reduced space $\widetilde{\mathcal{Y}}$.}
\end{algorithm}

\begin{remark}[Comparison with POD-greedy]
Note that in the POD-greedy algorithm \cite{HaaOhl08} one would consider the orthogonal projection on the reduced space $\widetilde{\mathcal{Y}}$ instead of the actual reduced solutions in \cref{eq:strongGreedyLambda}. However, for problems where the Galerkin projection deviates significantly from the orthogonal projection, we would expect that using the reduced solution gives superior results than the POD-greedy as the latter does not take into account the error due to the Galerkin projection which can be significant for instance close to resonances in a Helmholtz problem. We have performed numerical experiments for the same benchmark problem (parametrized Helmholtz equation) we consider in \cref{sec:Numerics} that confirm this conjecture.
\end{remark}

\subsection{Computational aspects of the fast-to-evaluate error estimator}\label{sec:online complexity}

At a first glance the complexity for evaluating $\mu\mapsto  \widetilde \Delta(\mu)$ is dominated by the solution of the $K$ reduced problems \cref{eq:monolithic}, meaning $K$ times the solution of a (dense) linear system of equations of size $n_{\widetilde{\mathcal{Y}}}$. The next proposition, inspired by Lemma 2.7 in \cite{zahm2017projection}, shows that one can actually evaluate $\mu\mapsto\widetilde\Delta(\mu)$ by solving only one linear system of size $n_{\widetilde{\mathcal{Y}}}$, which reduces the previous complexity by a factor $K$; the proof is provided in \cref{proof:monolithicAlternative}. 
Note however that the complexity for evaluating $\mu\mapsto\widetilde\Delta(\mu)$ is not completely independent on $K$. Indeed, as we employ the same reduced space for the approximation of $K$ dual problems, the dimension of $n_{\widetilde{\mathcal{Y}}}$ depends on $K$. The rate of the increase of $n_{\widetilde{\mathcal{Y}}}$ for growing $K$ will be investigated in numerical experiments in \cref{sec:Numerics}.

\begin{proposition}\label{prop:monolithicAlternative}
 The error indicator $\widetilde\Delta(\mu)$ defined by \eqref{eq: def a post est online} can be written as
 \begin{equation}\label{eq:rewrite_err}
  \widetilde\Delta(\mu) 
  = \left( \frac{1}{K}\sum_{i=1}^K  \big(Z_i^T \, \widetilde e(\mu) \big)^2 \right)^{1/2} ,
 \end{equation}
 where $ \widetilde e(\mu) \in\widetilde{\mathcal{Y}}$ is the solution to
 \begin{equation}\label{eq:monolithicAlternative}
   \widetilde e(\mu)\in\widetilde{\mathcal{Y}} \,, \quad
  \langle A(\mu) \widetilde e(\mu),v \rangle 
  = \langle r(\mu), v \rangle 
  \,, \quad \forall v\in\widetilde{\mathcal{Y}}.
 \end{equation}
 
\end{proposition}

Besides giving an alternative way of computing $\widetilde\Delta(\mu)$, \cref{prop:monolithicAlternative} also gives a new insight into the fast-to-evaluate error estimator. Reformulating Problem \cref{eq:monolithicAlternative} as
\begin{equation*}
 \widetilde e(\mu)\in\widetilde{\mathcal{Y}} \,, \quad
 \langle A(\mu) \big(\widetilde u(\mu) + \widetilde e(\mu) \big),v \rangle = \langle f(\mu), v \rangle 
 \,, \quad \forall v\in\widetilde{\mathcal{Y}} ,
\end{equation*}
demonstrates that $\widetilde e(\mu)\in\widetilde{\mathcal{Y}}$ may be interpreted as a correction of the primal approximation $\widetilde u(\mu)$, so that $\widetilde u(\mu) + \widetilde e(\mu)$ is an \textit{enriched solution} of the original problem \cref{eq:AUB} compared to $\widetilde u(\mu)$. Since $\widetilde{\mathcal{Y}}$ is not designed for improving the primal approximation $\widetilde u(\mu)$, one cannot reasonably hope that the correction $\widetilde e(\mu)$ improves significantly $\widetilde u(\mu)$. 
However the norm of $\widetilde e(\mu)$, estimated by the fast-to-evaluate error estimator \cref{eq:rewrite_err}, gives relevant information about the error $\|u(\mu)-\widetilde u(\mu)\|_\Sigma$. \blue{Finally, we emphasize again that the primal reduced space is in general not a subspace of the dual reduced space and that the intersection of the primal and dual reduced space can even be empty. As a consequence, the right-hand side in \eqref{eq:monolithicAlternative} is in general not zero for test functions from the dual reduced space.}

\begin{remark}\label{rmk:computationComplexityDelta}
 Assume $A(\mu) = \sum_{q=1}^{Q_{A}} \alpha_{q}(\mu) A_{q}$ and $f(\mu) = \sum_{q=1}^{Q_{f}} \zeta_{q}(\mu) f_{q}$ with $A_{q}, f_{q}$ parameter-independent and consider $\widetilde u(\mu)$ as the Galerkin projection onto some primal reduced order space $\widetilde{\mathcal{X}}$ of dimension $n_{\widetilde{\mathcal{X}}}$.
Since all inner products involving high dimensional quantities can be preassembled, the marginal computational complexity of $\widetilde\Delta(\mu)$ is $\mathcal{O}(Q_{A}n_{\widetilde{\mathcal{Y}}}^{2} + Q_{f}n_{\widetilde{\mathcal{Y}}} + Q_{A}n_{\widetilde{\mathcal{Y}}}n_{\widetilde{\mathcal{X}}})$ for assembling \cref{eq:monolithicAlternative}, $\mathcal{O}(n_{\widetilde{\mathcal{Y}}}^{3})$ for solving \cref{eq:monolithicAlternative} and $\mathcal{O}(K n_{\widetilde{\mathcal{Y}}})$ for calculating \cref{eq:rewrite_err}. 
 For moderate $Q_{A}$ the marginal computational complexity of $\widetilde\Delta(\mu)$ is thus dominated by $\mathcal{O}(n_{\widetilde{\mathcal{Y}}}^{3})$, \textit{i.e.} the costs for solving \cref{eq:monolithicAlternative}.
\end{remark}

\begin{remark}[Comparison with hierarchical type error estimators \cite{HORU18}]\label{rmk:hierarchicalEstimator}
 An alternative strategy for estimating the error is to measure the distance between the approximation $\widetilde u(\mu) \in \widetilde{\mathcal{X}}$ and a reference solution $u_{\mathrm{ref}}(\mu)$, which is an improved approximation of $u(\mu)$ compared to $\widetilde u(\mu)$. 
 When using projection based model order reduction $u_{\mathrm{ref}}(\mu)$ can be defined as a Galerkin projection onto an enriched reduced space of the form of $\widetilde{\mathcal{X}}+\widetilde{\mathcal{Y}}$, as proposed in \cite{HORU18}. Unlike our approach, the space $\widetilde{\mathcal{Y}}$ ought to be adapted for capturing the error $u(\mu) - \widetilde u(\mu)$. The complexity for evaluating such a hierarchical error estimator is dominated by the solution of a dense system of equations of size $\text{dim}(\widetilde{\mathcal{X}}+\widetilde{\mathcal{Y}})$. 
 In contrast, our approach requires the solution of a system of equations whose size is independent on the dimension of the primal reduced space $\widetilde{\mathcal{X}}$, see the above \cref{rmk:computationComplexityDelta}.
\end{remark}

\section{Numerical experiments}\label{sec:Numerics}

We numerically demonstrate various theoretical aspects of the proposed error estimator.
Our benchmark is a parameterized Helmholtz equation for which a reduced order solution is obtained by the RB method. 
Estimating the error in this reduced order model is challenging because, around the resonances, we lose the coercivity of the operator which makes a posteriori error estimation quite difficult with standard methods.

Let us mention here that all the training sets $\mathcal{P}^{\text{train}}$ (or $\mathcal{P}^{\text{train}}_K$) and all the online sets $\mathcal{S}$ are comprised of snapshots selected independently and uniformly at random in $\mathcal{P}$ (or in $\mathcal{P}_K$).
Those (random) sets are redrawn at each new simulation, unless mentioned otherwise.

\subsection{Benchmark: Multi-parametric Helmholtz equation}\label{subsect:Helmholtz}

\begin{figure}[t]
\begin{minipage}[c]{0.45\textwidth}
\centering
  \includegraphics[width=1\textwidth]{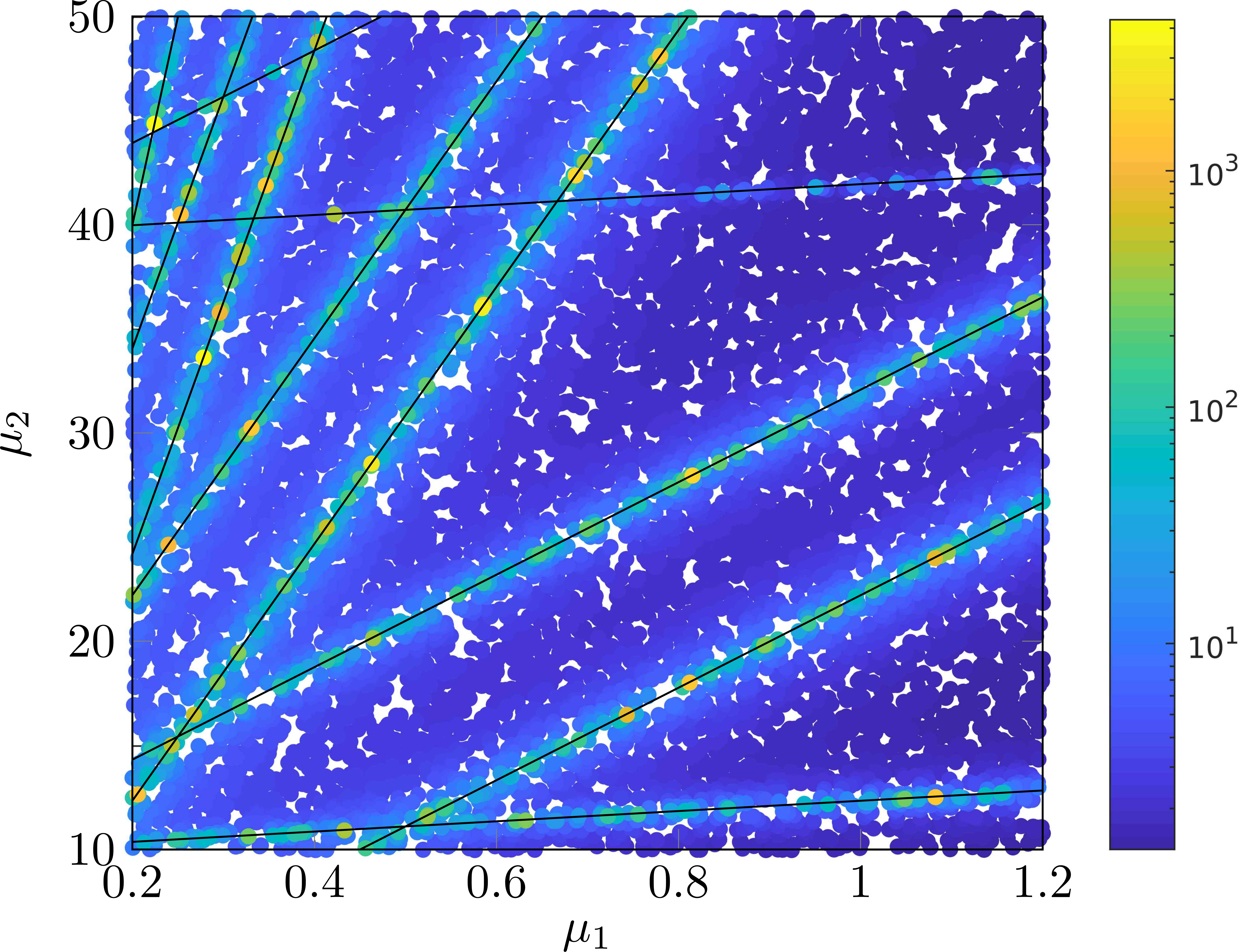}
\end{minipage}
\qquad
\begin{minipage}[c]{0.4\textwidth}
 \caption{Norm of the solution $\|u(\mu)\|_\Sigma$ \blue{over the online set $\mathcal{S}\subset\mathcal{P}$ with $\#\mathcal{S}=10^4$ and $\Sigma = R_X$}. The lines represent the resonances of the Helmholtz equation (computed analytically).} 
 \label{fig:ParamDomain}
\end{minipage}
 \vspace{-15pt}
\end{figure}

Consider the parameterized Helmholtz equation
\begin{alignat}{2}
\nonumber 
- \partial_{x_{1}x_{1}}\U - \mu_{1} \partial_{x_{2}x_{2}} \U - \mu_{2}\, \U &= \RHS & &\text{ in } D:=(0,1)\times (0,1),\\
\nonumber \U&=0 & &\text{ on } (0,1)\times \{0\},\\[-0.75em]
\label{eq:model problem} & & \\[-0.75em]
\nonumber \partial_{x_{2}} \U &= \cos(\pi x_1) \quad \quad & &\text{ on } (0,1)\times \{1\}, \\
\nonumber \partial_{x_{1}} \U &= 0 & &\text{ on } \{0,1\}\times (0,1).
\end{alignat}
The solution $\U = \U(\mu)$ is parameterized by $\mu=(\mu_{1},\mu_{2}) \in \mathcal{P}:=[0.2,1.2]\times [10,50]$, where $\mu_{1}$ accounts for anisotropy and $\mu_{2}$ is the wavenumber squared.
The source term $\RHS$ is defined by $\RHS(x_{1},x_{2})= \RHS_{1}(x_{1})\RHS_{2}(x_{2})$ for any $(x_1,x_2)\in D$, where
\begin{alignat*}{2}
\RHS_{1}(x_{1}) := \begin{cases} 
5 \quad &\text{ if } \enspace 0 \leq x_{1} \leq 0.1,\\
-5 \enspace &\text{ if } \enspace 0.2 \leq x_{1} \leq 0.3,\\
10 \enspace &\text{ if }\enspace 0.45 \leq x_{1} \leq 0.55,\\
-5 \enspace &\text{ if } \enspace 0.7 \leq x_{1} \leq 8,\\
5 \enspace &\text{ if } \enspace 0.9 \leq x_{1} \leq 1,\\
0 &\text{ else},
\end{cases} \quad \text{ and } 
\RHS_{2}(x_{2}) := \begin{cases}
1 \quad &\text{ if } \enspace 0.5 \leq x_{2} \leq 1, \\
0 &\text{ else}.
\end{cases}
\end{alignat*}
A similar test case with a smaller parameter set has been considered in \cite{Hetal10}.
The resonances can be determined analytically and are depicted by the black lines in \cref{fig:ParamDomain}. Because of the multi-parameter setting, we have resonance surfaces which are more difficult to deal with than a union of isolated resonance frequencies in the single-parameter setting; see \cite{Hetal10}. Moreover, we observe that in the region $[0.2,0.4] \times [30,50] \subset \mathcal{P}$ there are quite a few resonance surfaces that are also relatively close together, making this an even more challenging situation both for the construction of suitable reduced models and even more for a posteriori error estimation.  

We employ the Finite Element (FE) method to discretize the weak solution of \cref{eq:model problem}.
To that end, we define a FE space $X^{h} \subset X:=\{ \mathfrak{v} \in H^{1}(D)\, : \, \mathfrak{v}(x_{1},0) = 0\}$ by means of a regular mesh with square elements of edge length $h=0.01$ and FE basis functions that are piecewise linear in $x_{1}$ and $x_{2}$ direction, resulting in a FE space of $N = \dim(X^{h})=10100$. The FE approximation $\U^h(\mu)$ is defined as the Galerkin projection of $\U(\mu)$ on $X^h$, and we denote by $u(\mu)\in\mathbb{R}^{N}$ the vector containing the coefficients of $\U^h(\mu)$ when expressing it in the FE basis. Moreover, we denote by $R_X\in\mathbb{R}^{N\times N}$ the discrete Riesz map associated with the $H^1$-norm, which is such that $u(\mu) R_X u(\mu) = \| \U^h(\mu) \|_{H^{1}(D)}^2$ for any $\mu\in\mathcal{P}$. By default the covariance matrix $\Sigma$ is always chosen to be $\Sigma=R_X$, unless mentioned otherwise.

We may also consider a QoI defined as the trace of the FE solution on the boundary $\Gamma = \{0\}\times (0,1) \subset\partial D$, meaning $\U^h_{|\Gamma}(\mu)$.
We denote by $s(\mu) \in \mathbb{R}^{100}$ the vector containing those entries of $u(\mu) \in \mathbb{R}^{N}$ that are associated with the grid points on $\Gamma$. Then, we can write $s(\mu)=Lu(\mu)$ where $L\in\mathbb{R}^{100\times N}$ is an extraction matrix. To measure the error associated with the QoI, we use the norm $\|\cdot\|_W$ defined as the discretization of the $L^2(\Gamma)$ norm, which is such that $\|s(\mu)\|_{W} = \| \U_{|\Gamma}^h(\mu) \|_{L^2(\Gamma)}$ for any $\mu\in\mathcal{P}$.

The primal RB approximation $\widetilde u(\mu)$ is defined as the Galerkin projection of $u(\mu)$ onto the space of snapshots, meaning
$
 \widetilde u(\mu) \in \widetilde{\mathcal{X}} := \text{span}\{ u(\mu_1),u(\mu_2),\hdots \} ,
$
where the parameters $\mu_1,\mu_2,\hdots$ are selected in a greedy way based on the dual norm of the residual associated with \cref{eq:model problem}. Each time we run \cref{algo:qoi_greedy}, we use a reference solution $u_{\mathrm{ref}}(\mu)$ defined as an RB approximation of $u(\mu)$ using $n_{\widetilde{\mathcal{X}}} + 10$ basis functions, where $n_{\widetilde{\mathcal{X}}}:=\dim(\widetilde{\mathcal{X}})$.
Note that this reference solution appears only in the offline stage.

\subsection{Randomized a posteriori error estimation with exact dual}\label{subsec:numerics_exact_dual}

\begin{figure}[t]
    \centering
    \includegraphics[width=0.95\textwidth]{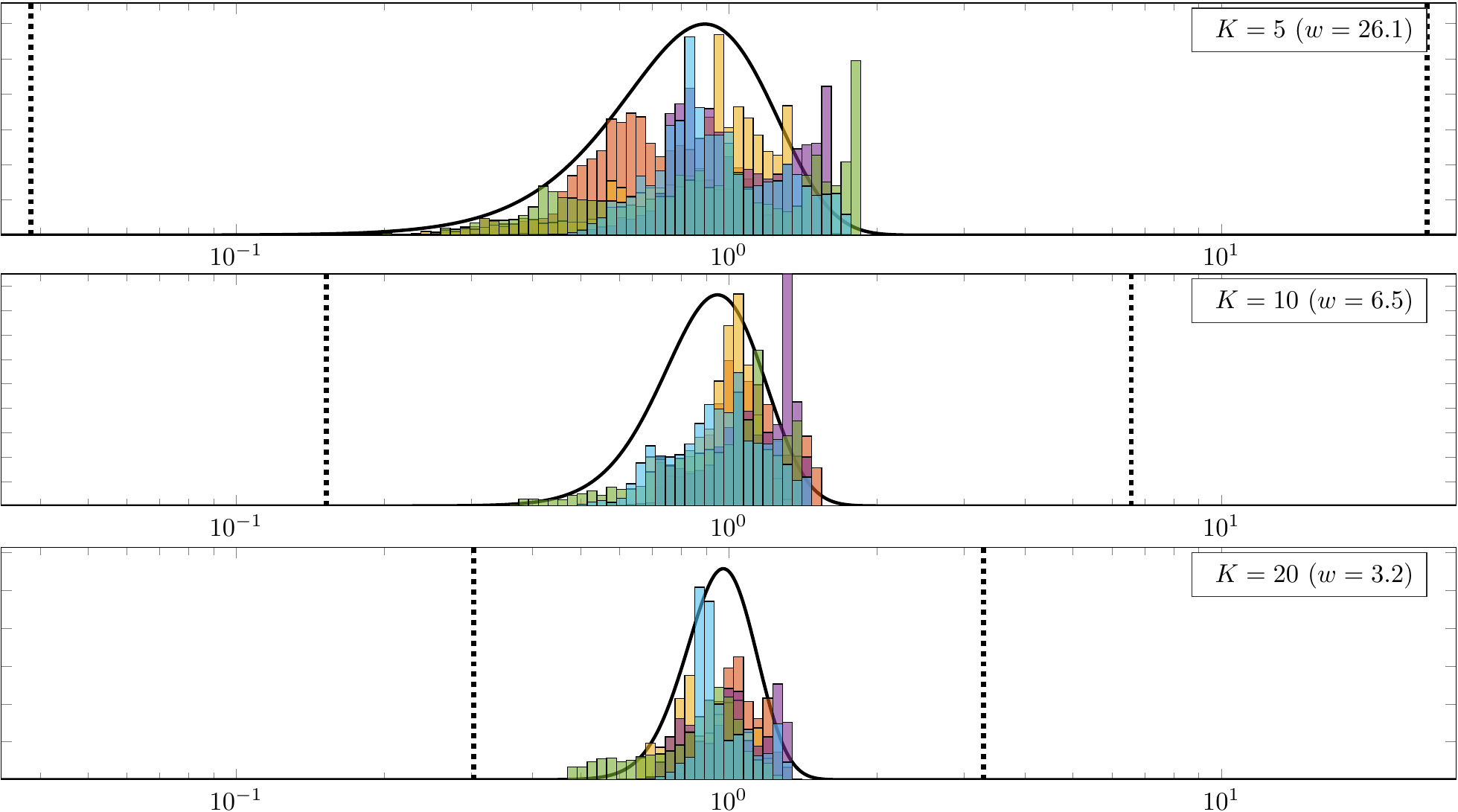}
  \caption{Histograms of $\{ \Delta(\mu)/\|u(\mu) - \widetilde u(\mu) \|_{\Sigma} \,,\,\mu\in\mathcal{S}\}$ for $n_{\widetilde{\mathcal{X}}}=10$ for five different realizations of the vectors $Z_1,\hdots,Z_K$, one color per realization. 
  Dashed lines: value of $1/w$ and $w$, where $w$ is obtained from \cref{coro:truth est Condition S} prescribing $\delta=10^{-2}$.}
  \label{fig:pdfXi2_Vs_DeltaExactDual}
  \vspace{-15pt}
\end{figure}

We demonstrate here the statistical properties of the error estimator $\Delta(\mu)$ defined by \cref{eq:DualTrick}.
\cref{fig:pdfXi2_Vs_DeltaExactDual} shows histograms of the effectivity indices $\{ \Delta(\mu)/\|u(\mu) - \widetilde u(\mu) \|_{\Sigma} \,,\,\mu\in\mathcal{S}\}$ for five different realizations of the vectors $Z_{1},\hdots,Z_{K}$. Here, the same online set $\mathcal{S}$ with $\#\mathcal{S}=10^{4}$ is used.
We observe that for each of the five realizations, the effectivity indices $ \Delta(\mu)/\|u(\mu) - \widetilde u(\mu) \|_{\Sigma}$ lie in the interval $[1/w , w]$ for any $\mu\in\mathcal{S}$, as predicted by \cref{coro:truth est S}.
This theoretical bound looks however pessimistic, as the effectivities for $K=5$ (resp. $K=10$) lie in the interval $[1/w,w]$ that corresponds to $K=10$ (resp. $K=20$).
This might be due to the rather crude union bound argument.

The solid lines on \cref{fig:pdfXi2_Vs_DeltaExactDual} represent the probability density function (pdf) of $\sqrt{Q/K}$ where $Q\sim \chi^2(K)$.
This is the pdf of $\Delta(\mu)/\|u(\mu) - \widetilde u(\mu) \|_{\Sigma}$ for any fixed $\mu$.
Even though the histograms depicted on \cref{fig:pdfXi2_Vs_DeltaExactDual} are not representing that pdf (instead they represent the distribution of the effectivity index among the set $\mathcal{S}$), we observe good accordance with the black line. In particular we observe a concentration phenomenon of the histograms around 1 when $K$ increases.

\subsection{Approximation of the dual problems}\label{subsec:numerics_approx_dual}
\subsubsection{Construction of the dual space}

\begin{figure}[h!] \centering
 \begin{subfigure}[c]{0.300\textwidth}\centering
  \includegraphics[width=\textwidth]{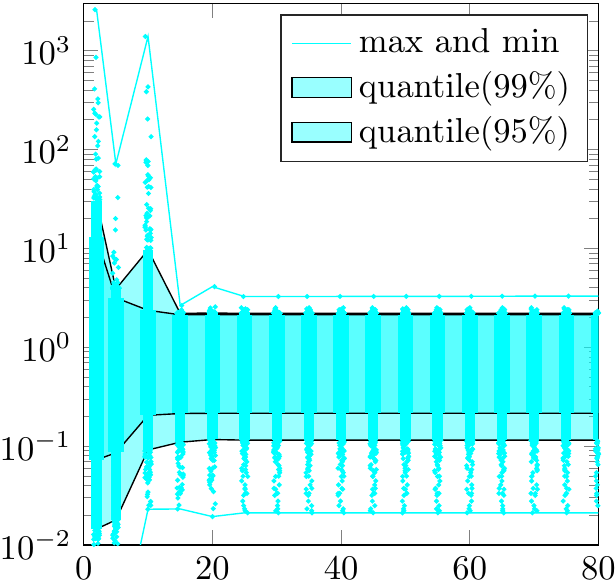}\vspace{-0.2cm}
  \caption{Alg.~\ref{algo:qoi_greedy}, $R_X$, $K=2$}\label{subfig:k2_alt}
 \end{subfigure} ~
 \begin{subfigure}[c]{0.300\textwidth}\centering
  \includegraphics[width=\textwidth]{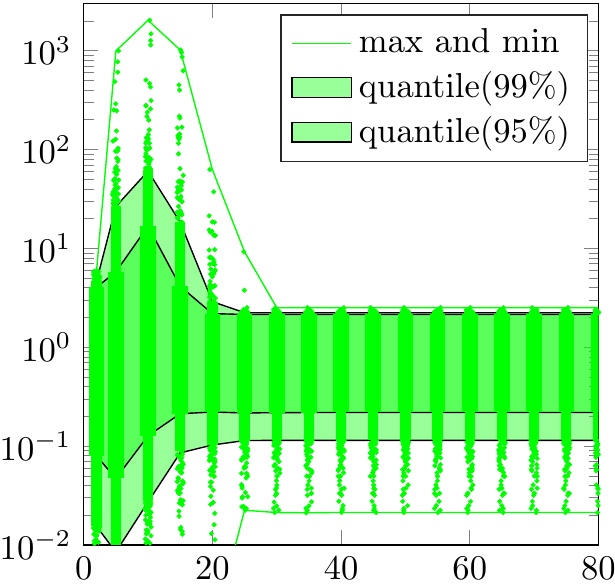}\vspace{-0.2cm}
  \caption{Alg.~\ref{algo:greedy}, $R_X$, $K=2$}
 \end{subfigure}
 \begin{subfigure}[c]{0.300\textwidth}\centering
  \includegraphics[width=\textwidth]{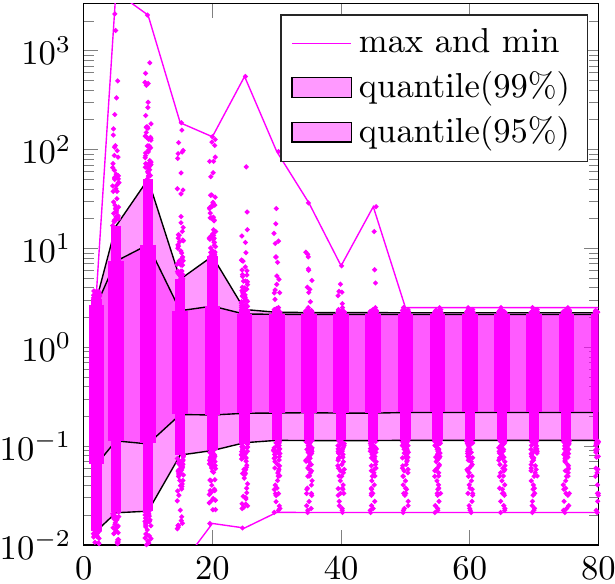}\vspace{-0.2cm}
  \caption{POD, $R_X$, $K=2$}
 \end{subfigure} 
 
 \vspace{-0.3cm}
 
 \begin{subfigure}[c]{0.300\textwidth}\centering
  \includegraphics[width=\textwidth]{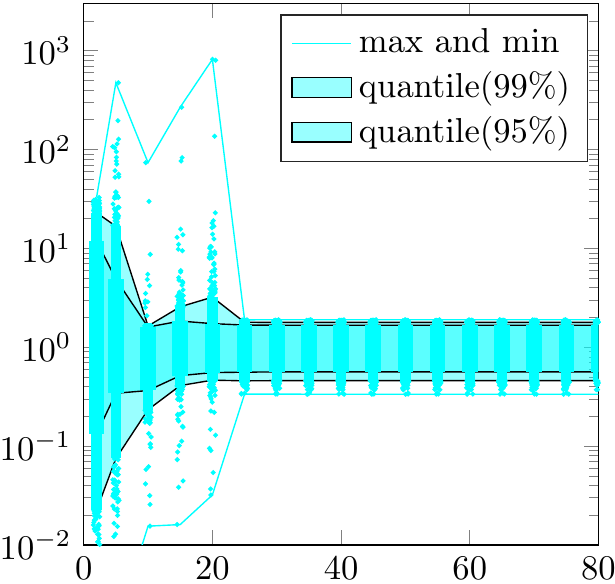}\vspace{-0.2cm}
  \caption{Alg.~\ref{algo:qoi_greedy}, $R_X$, $K=5$}\label{subfig:k5_alt}
 \end{subfigure} ~
 \begin{subfigure}[c]{0.300\textwidth}\centering
  \includegraphics[width=\textwidth]{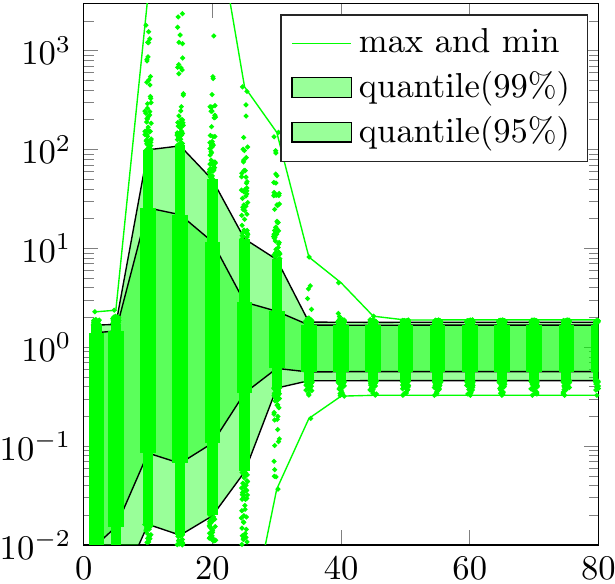}\vspace{-0.2cm}
  \caption{Alg.~\ref{algo:greedy}, $R_X$, $K=5$}\label{subfig:k5_stand}
 \end{subfigure}
 \begin{subfigure}[c]{0.300\textwidth}\centering
  \includegraphics[width=\textwidth]{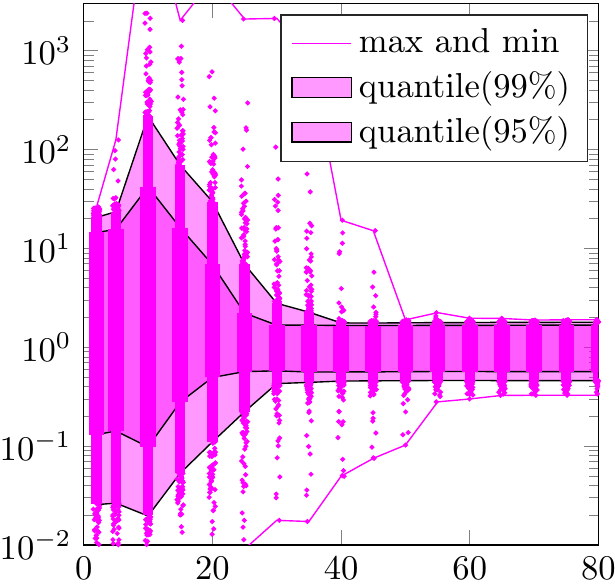}\vspace{-0.2cm}
  \caption{POD, $R_X$, $K=5$}\label{subfig:k5_pod}
 \end{subfigure} 
 
 \vspace{-0.3cm}
 
 \begin{subfigure}[c]{0.300\textwidth}\centering
  \includegraphics[width=\textwidth]{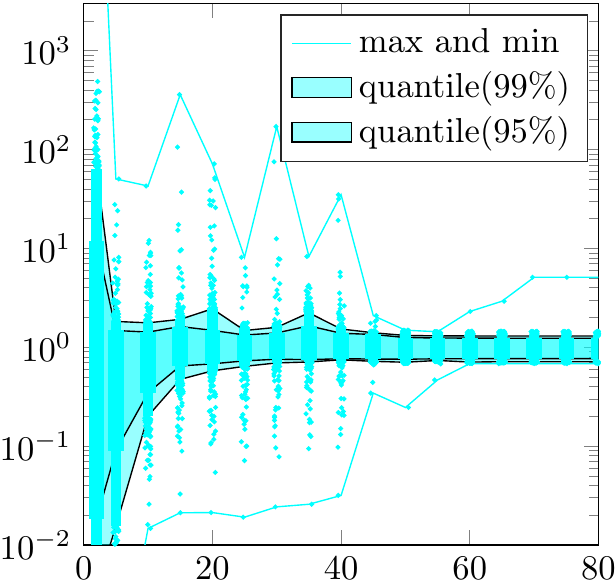}\vspace{-0.2cm}
  \caption{Alg.~\ref{algo:qoi_greedy}, $R_X$, $K=20$}
 \end{subfigure} ~
 \begin{subfigure}[c]{0.300\textwidth}\centering
  \includegraphics[width=\textwidth]{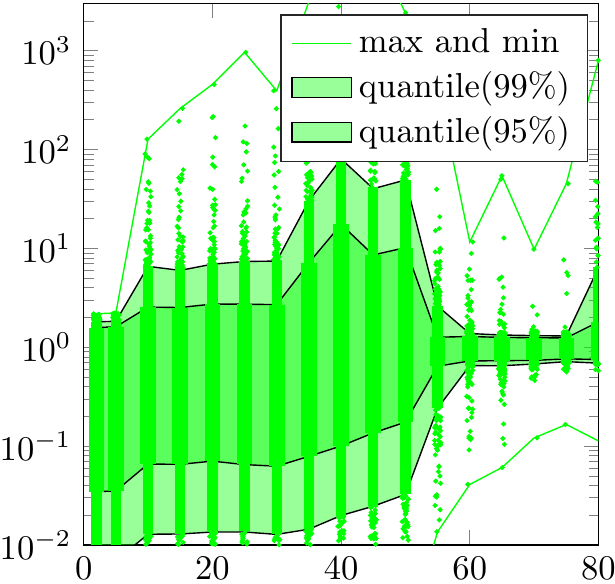}\vspace{-0.2cm}
  \caption{Alg.~\ref{algo:greedy}, $R_X$, $K=20$}
 \end{subfigure}
 \begin{subfigure}[c]{0.300\textwidth}\centering
  \includegraphics[width=\textwidth]{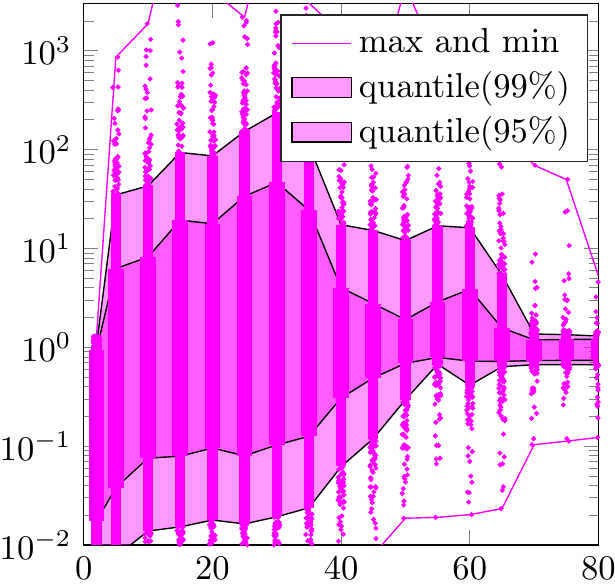}\vspace{-0.2cm}
  \caption{POD, $R_X$, $K=20$}\label{subfig:k20_pod}
  \end{subfigure} 
  
 \vspace{-0.3cm}
 
  \begin{subfigure}[c]{0.300\textwidth}\centering
  \includegraphics[width=\textwidth]{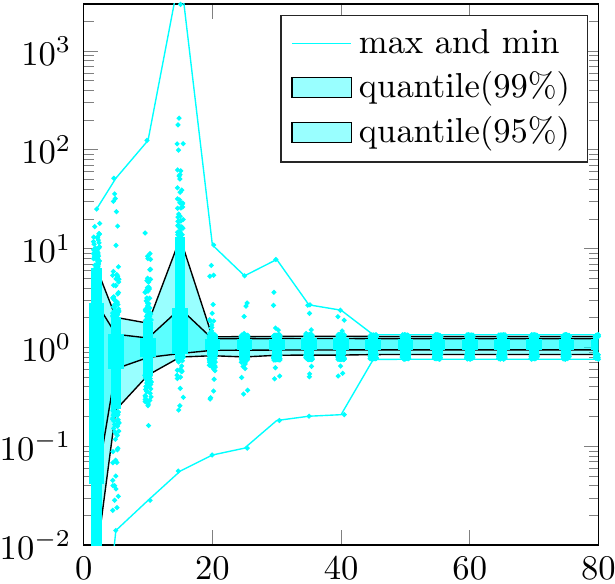}\vspace{-0.2cm}
  \caption{Alg.~\ref{algo:qoi_greedy}, QoI, $K=20$}\label{subfig:qoi_alt}
 \end{subfigure} ~
 \begin{subfigure}[c]{0.300\textwidth}\centering
  \includegraphics[width=\textwidth]{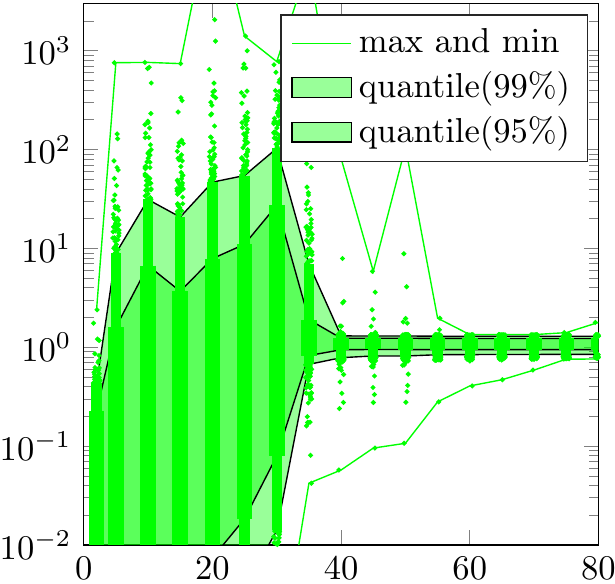}\vspace{-0.2cm}
  \caption{Alg.~\ref{algo:greedy}, QoI, $K=20$}\label{subfig:qoi_stand}
 \end{subfigure}
 \begin{subfigure}[c]{0.300\textwidth}\centering
  \includegraphics[width=\textwidth]{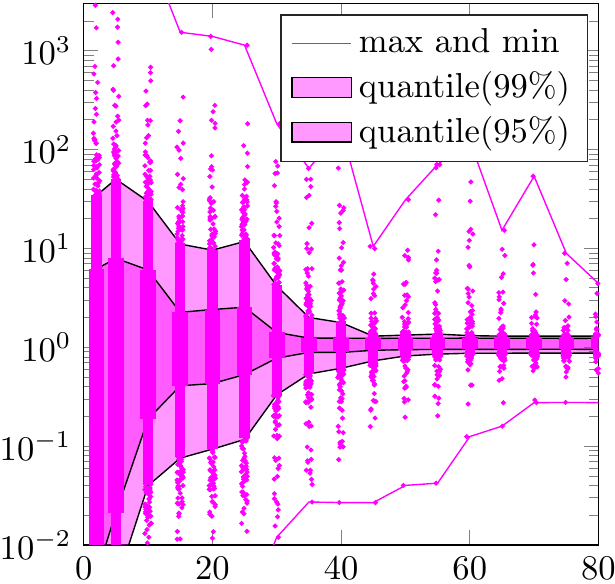}\vspace{-0.2cm}
  \caption{POD, QoI, $K=20$}\label{subfig:qoi_pod}
 \end{subfigure} 
 
 \vspace{-0.6cm}
 
 \caption{Maximum, minimum, and two quantiles ($99\%$ and $95\%$) of $\{ \widetilde\Delta(\mu)/\|u(\mu)-\widetilde u(\mu)\|_{\Sigma} : \mu\in\mathcal{S}\}$ as a function of $n_{\widetilde{\mathcal{Y}}}$.
 The dual reduced space $\widetilde{\mathcal{Y}}$ is constructed by \cref{algo:qoi_greedy} (left column), \cref{algo:greedy} (middle column), and POD (right column).
 The first three rows corresponds to different values of $K=2,5,20$ with $\Sigma=R_X$.
 The last row corresponds to $\Sigma = L^TR_WL$ (the QoI) with $K=20$.
 On each row we use the same realization of the vectors $Z_1,\hdots,Z_K$, which allows a fair comparison of the different algorithms.
 For each plot we use $\#\mathcal{P}^\text{train}=10^3$ and $\#\mathcal{S}=10^4$.
 }
 \label{fig:Approx_Dual_RX}
\end{figure}

In \cref{fig:Approx_Dual_RX} we compare the maximum, the minimum, the 95\% quantile and the 99\% quantile of $\{ \widetilde\Delta(\mu)/\|u(\mu)-\widetilde u(\mu)\|_{\Sigma} : \mu\in\mathcal{S}\}$ where the dual reduced space is constructed either by \cref{algo:greedy} (with $\|\cdot\|_{*}=\|\cdot\|_{R_X^{-1}}$), by \cref{algo:qoi_greedy} or by a POD.
We observe that by using \cref{algo:qoi_greedy} we need many fewer dual basis functions than for \cref{algo:greedy} and for the POD. 
In detail, we see for instance in \cref{subfig:k5_alt}, \cref{subfig:k5_stand} and \cref{subfig:k5_pod} that for $K=5$ employing \cref{algo:qoi_greedy} requires about $n_{\widetilde{\mathcal{Y}}}=20$ dual basis functions to have 99\% of the samples in the interval $[1/3,3]$, while when using the \cref{algo:greedy} or the POD we need about $35$ or $30$ basis functions, respectively. We emphasize that for $K=20$ the difference is even larger. 
\blue{Moreover,} when considering the QoI (last row of \cref{fig:Approx_Dual_RX}), the difference between \cref{algo:qoi_greedy} and the \cref{algo:greedy} and POD is less pronounced but still considerable. This significant disparity can be explained by the fact that while both the POD and the \cref{algo:greedy} try to approximate the $K$ dual solutions $\widetilde Y_{1}(\mu),\hdots,\widetilde Y_{K}(\mu)$, \cref{algo:qoi_greedy} is driven by the approximation of the error estimator $\Delta(\mu)$ and thus a scalar quantity; compare the selection criteria \cref{eq:weakGreedy_training} and \cref{eq:strongGreedy}. 
This also explains why the discrepancy increases significantly for growing $K$: While POD and \cref{algo:greedy} have to approximate a more complex object (the $K$ dual solutions), we only obtain an additional summand in $\widetilde \Delta(\mu)$ for each additional random right-hand side. 
Let us also highlight the significant difference between the maximum value and the $99\%$ quantile over the parameter set and the somewhat erratic behavior of the maximum, which both seem to be due to the resonance surfaces. As indicated above this motivates considering for instance the $99\%$ quantile as a stopping criterion in both \cref{algo:greedy} and \cref{algo:qoi_greedy}.

\subsubsection{Dimension of the dual space}

\begin{figure}[t]
 \centering
 \begin{subfigure}[c]{0.3\textwidth}\centering
  \includegraphics[width=\textwidth]{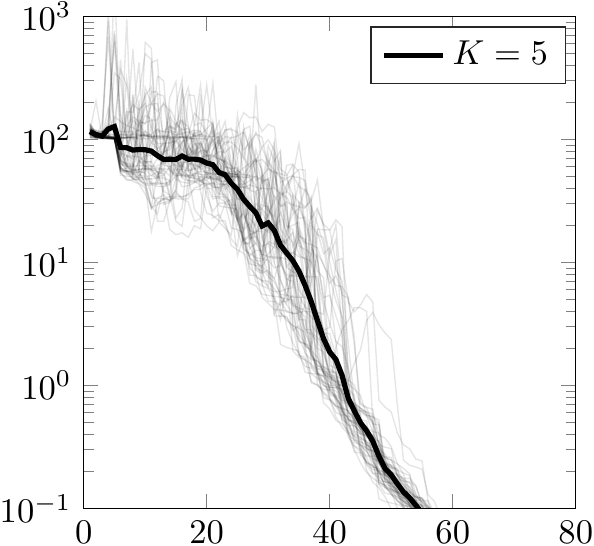}
 \end{subfigure}
 \begin{subfigure}[c]{0.3\textwidth}\centering
  \includegraphics[width=\textwidth]{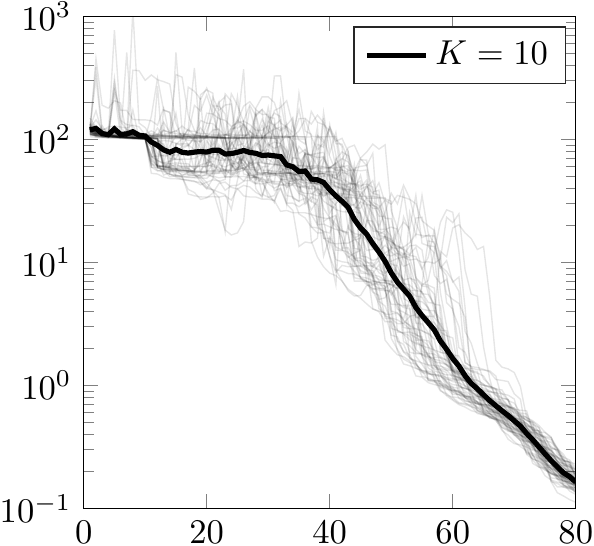}
 \end{subfigure}
 \begin{subfigure}[c]{0.3\textwidth}\centering
  \includegraphics[width=\textwidth]{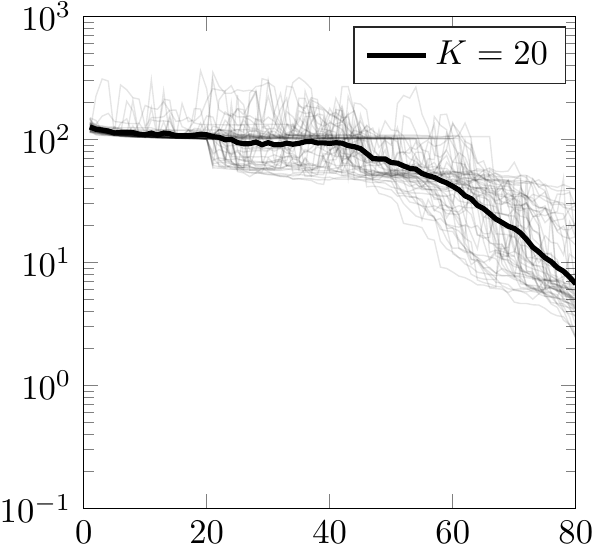}
 \end{subfigure} 
 \caption{Evolution of the $97.5\%$-quantile of $\{ \| A(\mu)^T \widetilde Y_i(\mu) - Z_i \|_*  :  (i,\mu)\in\mathcal{P}_K^\text{train} \}$ during the first 80 iterations of \cref{algo:greedy}. Each grey line corresponds to one realization of $Z_1,\hdots,Z_K$ and the black lines are the mean of the grey lines.}
 \label{fig:ConvGreedyDual_standard}
 \vspace{-15pt}
\end{figure}

\cref{tab:stand_greedy} shows statistics of the dimension of the dual reduced space $\widetilde{\mathcal{Y}}$ obtained by \cref{algo:greedy} with different stopping criterion. \blue{We consider $\tol=0.5$ and relax this tolerance by multiplication with a varying relaxation parameter $\rho$ taking the values $1,10$ and $100$. We observe that except for $\tol\cdot\rho=50$ and moderate $K$ the dimension of the dual space is in general quite large. Comparing with \cref{fig:Approx_Dual_RX}, we observe that choosing $\tol\cdot\rho=50$ is albeit sufficient to obtain an effectivity close to $1$.
Notice however that the use of \cref{cor:dualErrorAdditive} requires $\tol\cdot\rho \approx \varepsilon \leq 1/w \leq 1$, which excludes $\tol\cdot\rho=5$ and $\tol\cdot\rho=50$.} \cref{fig:ConvGreedyDual_standard} shows the evolution of the stopping criteria during the first 80 iterations of \cref{algo:greedy}. 
We observe a significant impact of $K$ on the convergence profiles: with $K=20$ the curves do not attain the tolerance $\tol=0.5$, which explains the results we observed in \cref{tab:stand_greedy}.

\blue{
\begin{table}
\centering
\begin{subtable}{.9\textwidth}
\centering
\footnotesize
\begin{tabular}{ |c|c | c | c| }
 \cline{2-4}
 \multicolumn{1}{c|}{~} & \blue{$\tol\cdot\rho=0.5$ ($\rho=1$)} & \blue{$\tol\cdot\rho=5$ ($\rho=10$)} & \blue{$\tol\cdot\rho=50$ ($\rho=100$)} \\ \hline
$K=5$ & $47.7 \,(\pm 1.9)$ & $39.9 \,(\pm 3.24)$ & $30.6 \,(\pm 9.3)$ \\
$K=10$ & $74.2 \,(\pm 1.75)$ & $56.2 \,(\pm 4.4)$ & $43.5 \,(\pm 10.8)$ \\
$K=20$ & $>80$ & $>80$ & $65.4 \,(\pm 9.65)$ \\
$K=50$ & $>80$ & $>80$ & $>80$ \\
 \hline
\end{tabular}
\caption{Stopping criterion: max ($q=1$) over $\mathcal{P}^\text{train}$ is \blue{$\leq \tol\cdot\rho$} }
\end{subtable}
\begin{subtable}{.9\textwidth}
\centering
\footnotesize
\begin{tabular}{ |c|c | c | c| }
 \cline{2-4}
 \multicolumn{1}{c|}{~} & \blue{$\tol\cdot\rho=0.5$ ($\rho=1$)} & \blue{$\tol\cdot\rho=5$ ($\rho=10$)} & \blue{$\tol\cdot\rho=50$ ($\rho=100$)} \\ \hline
$K=5$ & $45.4 \,(\pm 2.14)$ & $35.9 \,(\pm 3.29)$ & $16.9 \,(\pm 6.99)$ \\
$K=10$ & $70.6 \,(\pm 1.83)$ & $53.7 \,(\pm 4.17)$ & $31 \,(\pm 9.4)$ \\
$K=20$ & $>80$ & $>80$ & $55.1 \,(\pm 9.43)$ \\
$K=50$ & $>80$ & $>80$ & $>80$ \\
 \hline
\end{tabular}
\caption{Stopping criterion: $97.5\%$-quantile ($q=0.975$) over $\mathcal{P}^\text{train}$ is \blue{$\leq \tol\cdot\rho$}}
\end{subtable}
\caption{Mean ($\pm$standard deviation) of $n_{\widetilde{\mathcal{Y}}}$ over $100$ realizations of the $K$ vectors $Z_1,\hdots,Z_K$. Here $\widetilde{\mathcal{Y}}$ is built using \cref{algo:greedy} with different stopping criterion, \textit{i.e.} with different values for $q$ \blue{and $\rho$.}}\label{tab:stand_greedy}
\end{table}}

\begin{table}
\centering
\begin{subtable}{.9\textwidth}
\centering
\footnotesize
\begin{tabular}{ |c|c | c | c|c| }
 \cline{2-5}
 \multicolumn{1}{c|}{~} & $\tol=1.5$ & $\tol=2$ &$\tol=3$ & $\tol=5$ \\ \hline
$K=5$ & $26.9 \,(\pm 4.41)$ & $23 \,(\pm 5.58)$ & $19 \,(\pm 6.35)$ & $14.8 \,(\pm 6.52)$ \\
$K=10$ & $28.1 \,(\pm 7.64)$ & $22.1 \,(\pm 5.72)$ & $16.4 \,(\pm 5.25)$ & $12.3 \,(\pm 3.56)$ \\
$K=20$ & $28.4 \,(\pm 10.2)$ & $22.4 \,(\pm 9.56)$ & $16.3 \,(\pm 8.38)$ & $13.1 \,(\pm 8.04)$ \\
$K=50$ & $31.8 \,(\pm 11.3)$ & $21.8 \,(\pm 8.06)$ & $14.9 \,(\pm 4.93)$ & $11.3 \,(\pm 2.68)$ \\
\hline
\end{tabular}
\caption{Stopping criterion: max of $\alpha$ ($q=1$) over $\mathcal{P}^\text{train}$ is $\leq tol$, $ n_{\widetilde{\mathcal{X}}}=10$}
\end{subtable}

\begin{subtable}{.9\textwidth}
\centering
\footnotesize
\begin{tabular}{ |c|c | c | c|c| }
 \cline{2-5}
 \multicolumn{1}{c|}{~} & $\tol=1.5$ & $\tol=2$ &$\tol=3$ & $\tol=5$ \\ \hline
$K=5$ & $18.7 \,(\pm 4.89)$ & $13.9 \,(\pm 4.23)$ & $9.7 \,(\pm 4.25)$ & $6.66 \,(\pm 3.22)$ \\
$K=10$ & $18.2 \,(\pm 5)$ & $12 \,(\pm 3.51)$ & $7.64 \,(\pm 2.08)$ & $6.08 \,(\pm 1.64)$ \\
$K=20$ & $21.9 \,(\pm 6.96)$ & $13.9 \,(\pm 4.18)$ & $8.88 \,(\pm 2.56)$ & $6.02 \,(\pm 1.9)$ \\
$K=50$ & $25.1 \,(\pm 9.77)$ & $15.7 \,(\pm 5.74)$ & $9.44 \,(\pm 3.44)$ & $6.12 \,(\pm 2.03)$ \\
\hline
\end{tabular}
\caption{Stopping criterion: $97.5\%$-quantile ($q=0.975$) over $\mathcal{P}^\text{train}$ is $\leq tol$, $ n_{\widetilde{\mathcal{X}}}=10$}
\end{subtable}

\begin{subtable}{.9\textwidth}
\centering
\footnotesize
\begin{tabular}{ |c|c | c | c|c| }
 \cline{2-5}
 \multicolumn{1}{c|}{~} & $q=100\%$ $(\max)$ & $q=99\%$ & $q=97.5\%$ & $q=95\%$ \\ \hline
$K=5$ & $19 \,(\pm 6.35)$ & $11.7 \,(\pm 5.14)$ & $9.7 \,(\pm 4.25)$ & $8.52 \,(\pm 3.84)$ \\
$K=10$ & $16.4 \,(\pm 5.25)$ & $9.5 \,(\pm 2.59)$ & $7.64 \,(\pm 2.08)$ & $6.8 \,(\pm 2.01)$ \\
$K=20$ & $16.3 \,(\pm 8.38)$ & $10.5 \,(\pm 3.29)$ & $8.88 \,(\pm 2.56)$ & $7.38 \,(\pm 2.55)$ \\
$K=50$ & $14.9 \,(\pm 4.93)$ & $10.7 \,(\pm 3.46)$ & $9.44 \,(\pm 3.44)$ & $7.16 \,(\pm 2.58)$ \\ \hline 
\end{tabular}
\caption{Stopping criterion: $q$-quantile over $\mathcal{P}^\text{train}$ is $\leq \tol=3$, $ n_{\widetilde{\mathcal{X}}}=10$}
\end{subtable}

\begin{subtable}{.9\textwidth}
\centering
\footnotesize
\begin{tabular}{ |c|c | c | c| }
 \cline{2-4}
 \multicolumn{1}{c|}{~} & $ n_{\widetilde{\mathcal{X}}}=10$ & $ n_{\widetilde{\mathcal{X}}}=20$ & $ n_{\widetilde{\mathcal{X}}}=30$ \\ \hline
$K=5$ & $9.7 \,(\pm 4.25)$ & $6.74 \,(\pm 1.65)$ & $7.56 \,(\pm 2.27)$ \\
$K=10$ & $7.64 \,(\pm 2.08)$ & $9.86 \,(\pm 2.55)$ & $9.62 \,(\pm 3.1)$ \\
$K=20$ & $8.88 \,(\pm 2.56)$ & $14.2 \,(\pm 3.57)$ & $13.6 \,(\pm 2.7)$ \\
$K=50$ & $9.44 \,(\pm 3.44)$ & $22.1 \,(\pm 5.22)$ & $23.1 \,(\pm 5.27)$ \\
\hline
\end{tabular}
\caption{Stopping criterion: $97.5\%$-quantile ($q=0.975$) over $\mathcal{P}^\text{train}$ is $\leq \tol = 3$}\label{subtable:algo 2 dimX}
\end{subtable}
\caption{Mean ($\pm$standard deviation) of $n_{\widetilde{\mathcal{Y}}}$ over $100$ realizations of the $K$ vectors $Z_1,\hdots,Z_K$. Here $\widetilde{\mathcal{Y}}$ is built using \cref{algo:qoi_greedy} with different stopping criterion $q$ and $\tol$, and with different primal approximation $n_{\widetilde{\mathcal{X}}}=10,20,30$. Here $\#\mathcal{P}^\text{train}=10^4$.} \label{table:algo 2}
\end{table}

\begin{figure}
 \centering
 \begin{subfigure}[c]{0.3\textwidth}\centering
  \includegraphics[width=\textwidth]{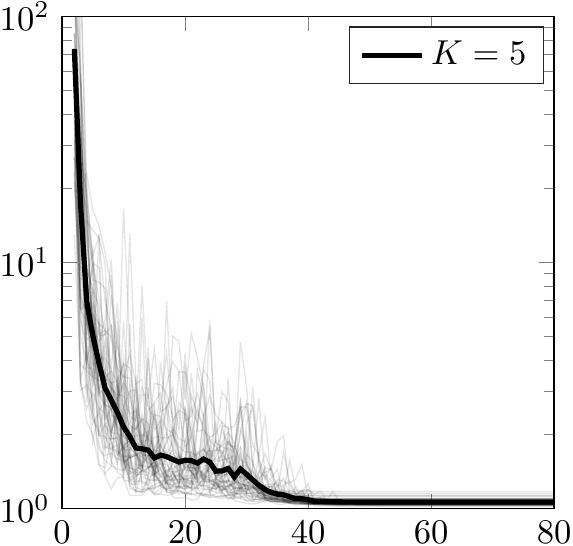}
 \end{subfigure}
 \begin{subfigure}[c]{0.3\textwidth}\centering
  \includegraphics[width=\textwidth]{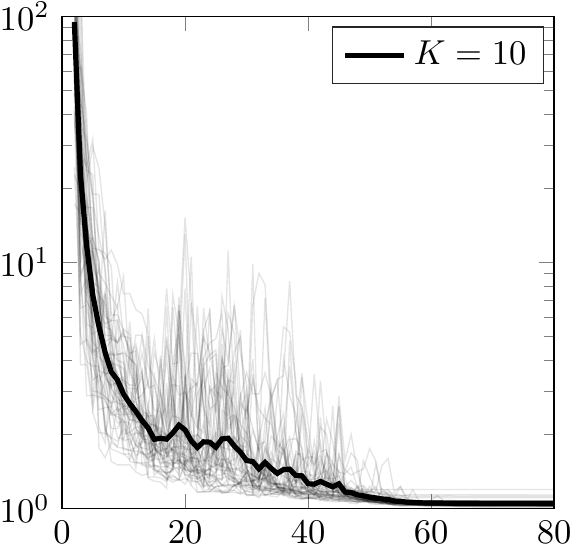}
 \end{subfigure}
 \begin{subfigure}[c]{0.3\textwidth}\centering
  \includegraphics[width=\textwidth]{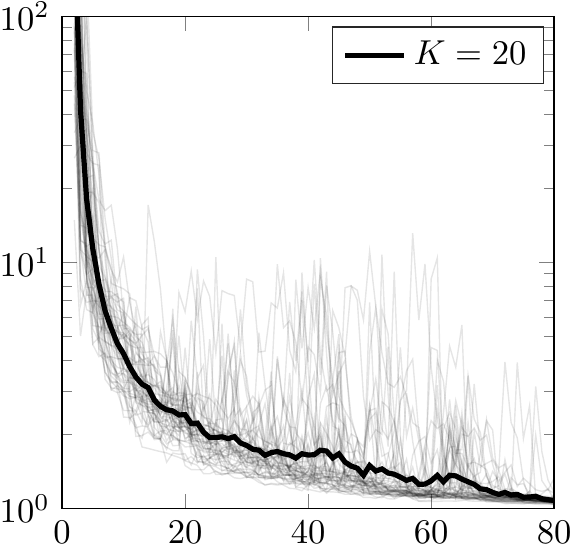}
 \end{subfigure} 
\caption{Evolution of the $97.5\%$-quantile of $\{ \max\{  \frac{\Delta_{\mathrm{ref}}(\mu)}{\widetilde\Delta(\mu)} ;\frac{ \widetilde\Delta(\mu)}{\Delta_{\mathrm{ref}}(\mu)}  \}  : \mu\in\mathcal{P}^\text{train} \}$ during the first 80 iterations of \cref{algo:greedy}. Here, $n_{\widetilde{\mathcal{X}}}=30$. Each grey line corresponds to one realization of $Z_1,\hdots,Z_K$ and the black lines are the mean of the grey lines.}
 \label{fig:ConvGreedyDual}
 \vspace{-15pt}
\end{figure}

In comparison, \cref{algo:qoi_greedy} yields much smaller dual reduced spaces; compare \cref{tab:stand_greedy} and \cref{table:algo 2}. We see in \cref{table:algo 2} that, except for $\tol=1.5$ the dimension of the dual RB space $\widetilde{\mathcal{Y}}$ is smaller than the dimension of the primal RB space $\widetilde{\mathcal{X}}$ when using the $95\%$, $97.5\%$, $99\%$-quantiles for the stopping criterion. Moreover, for instance for $\tol=3$ we see that for $n_{\widetilde{\mathcal{X}}}=20,30$ we can use (significantly) less dual than primal basis functions. However, we also see that tight tolerances for $\tol$ will lead in general to dual reduced spaces that have a larger dimension than the primal RB space. 
As larger tolerances $\geq 5$ may lead to an significant underestimation of the error (see \cref{fig:varying_alpha_w}), tolerances for $\tol$ between $1$ and $4$ seem to be preferable. \cref{fig:ConvGreedyDual} shows the evolution of the stopping criteria during the first 80 iterations of \cref{algo:qoi_greedy}. Note that for higher tolerances for $\tol$ it may happen for a realization that \cref{algo:qoi_greedy} terminates in a valley between two peaks.

Furthermore, we observe in \cref{table:algo 2} a very large standard deviation of about $10$ if we consider the maximum over the offline training set, while for the $95\%$,$97.5\%$,$99\%$ quantiles we have often a standard deviation of about $2$. Additionally, the dimension of the dual reduced spaces for the maximum is much larger than for the considered quantiles, but among the considered quantiles we observe only very moderate changes. Again, it seems that this behavior is due to the resonance surfaces. Moreover, as we obtain a very satisfactory effectivity of $\widetilde \Delta(\mu)$ when we use for instance the $99\%$ quantile (see \cref{subsub:numerics_online_perf}), we conclude that using quantiles between $97.5\%$ and $99\%$ as a stopping criterion in \cref{algo:qoi_greedy} seems advisable.

Finally, we observe both a very moderate dependency of the dimension of the dual reduced space constructed by \cref{algo:qoi_greedy} on $K$ and a rather mild dependency on $n_{\widetilde{\mathcal{X}}}$. Therefore, we conjecture that the proposed error estimator might also be applied rather complex problems.

\subsubsection{Performance of $\widetilde \Delta(\mu)$ on an online parameter set}\label{subsub:numerics_online_perf}

\begin{figure}
 \centering
 \includegraphics[width=0.95\textwidth]{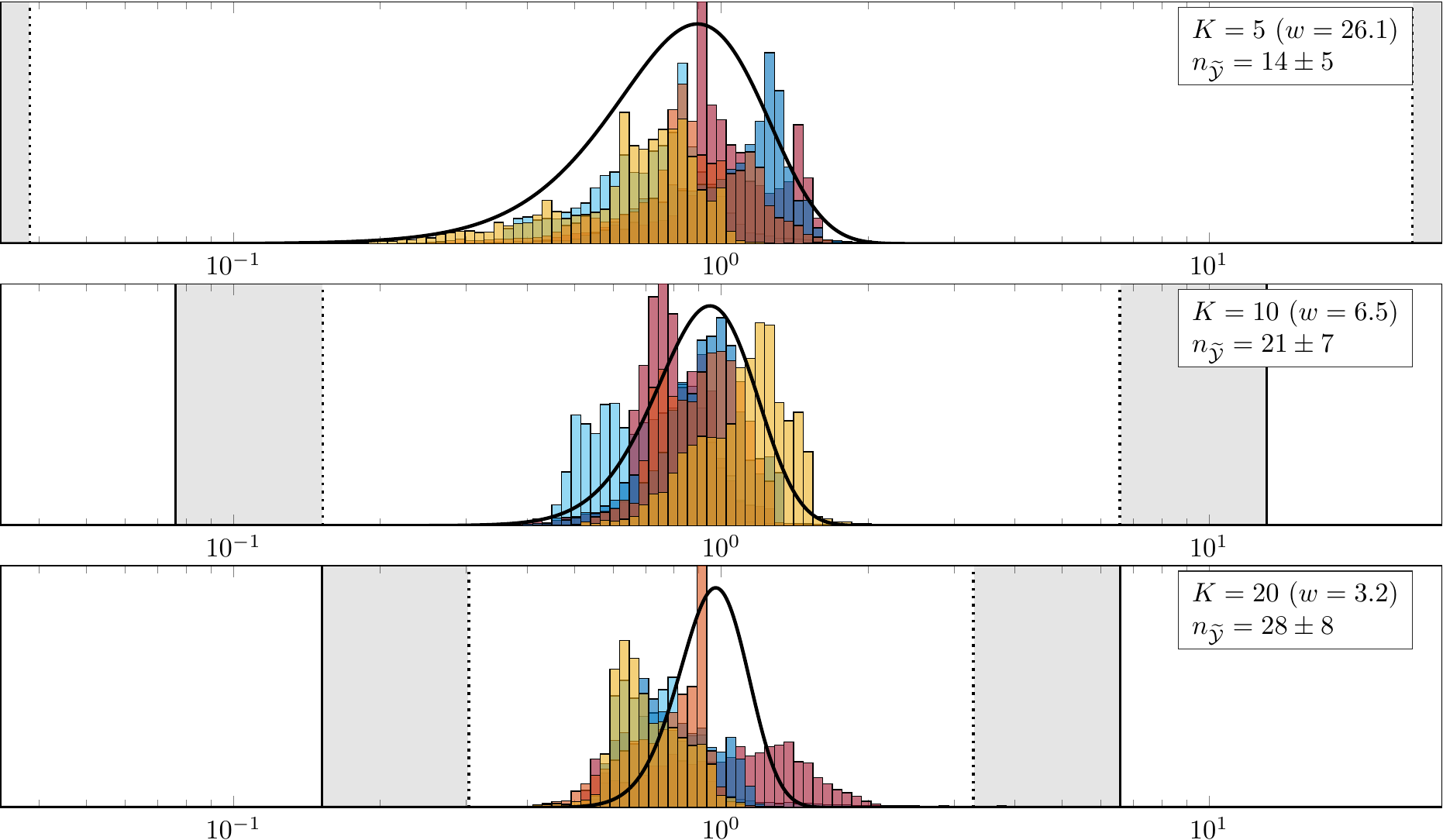}
 \caption{
 Histograms of $\{ \widetilde\Delta(\mu)/\|u(\mu)-\widetilde u(\mu)\|_{\Sigma} : \mu\in\mathcal{S}\}$ with $\#\mathcal{S}=10^4$ for $5$ realizations of $K$ random vectors $Z_{1},\hdots,Z_{K}$ (one color per realization). The dual reduced space $\widetilde{\mathcal{Y}}$ is built using \cref{algo:qoi_greedy} with $q=0.99$, $\tol=2$ and $\#\mathcal{P}^\text{train}=10^3$. Here $n_{\widetilde{\mathcal{X}}}=20$ and $\Sigma = R_X$. 
 The vertical dashed lines corresponds to $w^{-1}$ and $w$ where $w$ is obtained from \cref{coro:truth est Condition S} prescribing $\delta=10^{-2}$.
 The gray area corresponds to the amplification of the confidence interval due to $\alpha\approx\tol$, see \cref{prop:dualErrorMultiplicative}.}
 \label{fig:histograms_online}
 \vspace{-15pt}
\end{figure}

On \cref{fig:histograms_online} we plot the histograms of $\{ \widetilde\Delta(\mu)/\|u(\mu)-\widetilde u(\mu)\|_{\Sigma} : \mu\in\mathcal{S}\}$ for $5$ realizations of the random vectors $Z_{1},\hdots,Z_{K}$, where the dual reduced space $\widetilde{\mathcal{Y}}$ is built via \cref{algo:qoi_greedy}. 
We observe a similar behavior as for the error estimator $\Delta(\mu)$ with the exact dual, see \cref{fig:pdfXi2_Vs_DeltaExactDual}. In particular for all $\mu\in\mathcal{S}$ the effectivity index $\widetilde\Delta(\mu)/\|u(\mu)-\widetilde u(\mu)\|_{\Sigma}$ lies  between $(\alpha w)^{-1}$ and $(\alpha w)$, see \cref{prop:dualErrorMultiplicative}, where $\alpha$ is estimated by $\tol=2$.
Finally, we highlight that \cref{fig:histograms_online} demonstrates that with near certainty we obtain effectivities near unity with a dual space dimension on the same order as (or less than) the primal space dimension. Hence the costs for the a posteriori error estimator are about the same as those for constructing the primal approximation.

\begin{figure}
 \centering
 \includegraphics[width=0.95\textwidth]{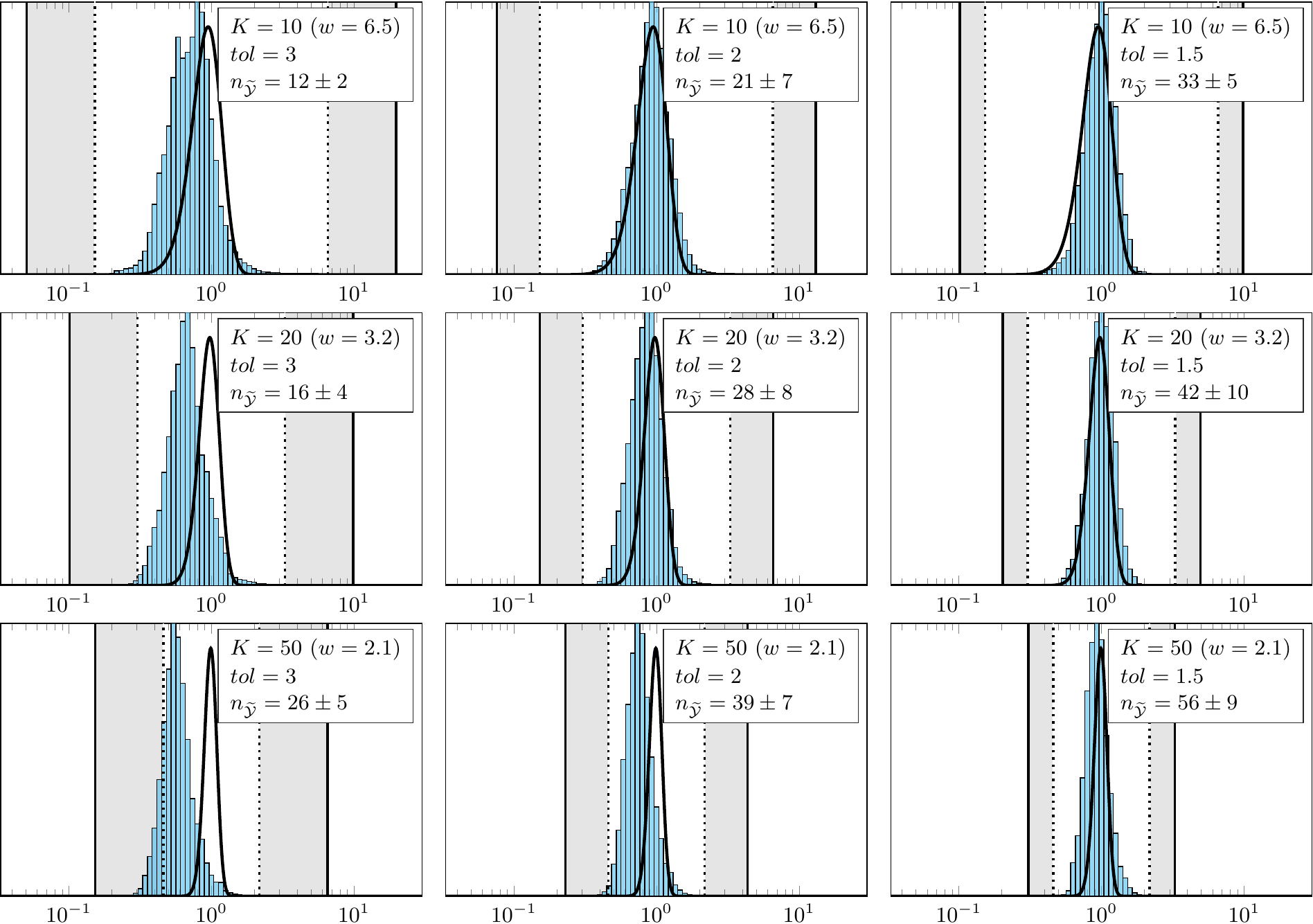}
 \caption{Histograms of the concatenation of 100 realizations of $\{ \widetilde\Delta(\mu)/\|u(\mu)-\widetilde u(\mu)\|_{\Sigma} : \mu\in\mathcal{S}\}$ where at each realization, the vectors $Z_1,\hdots,Z_K$, the training set $\mathcal{P}^\text{train}$ and the online set $\mathcal{S}$ are redrawn and $\widetilde{\mathcal{Y}}$ is rebuilt using \cref{algo:qoi_greedy} with $q=0.99$.
 The solid lines are the pdf of $\sqrt{Q/K}$ where $Q\sim \chi^2(K)$.
 As on \cref{fig:histograms_online}, the grey area corresponds to the amplification of the confidence interval $[(\alpha w)^{-1},(\alpha w)]$ due to $\alpha\approx\tol$.
 Here $n_{\widetilde{\mathcal{X}}}=20$, $\Sigma = R_X$, $\delta=10^{-2}$, $\#\mathcal{P}^\text{train}=10^3$ and $\#\mathcal{S}=10^4$.
 }
 \label{fig:varying_alpha_w}
 \vspace{-15pt}
\end{figure}

In order to understand the average performance of the online-efficient error indicator, we plot in \cref{fig:varying_alpha_w} the histograms of the concatenation of 100 realizations of the effectivity indices $\{ \widetilde\Delta(\mu)/\|u(\mu)-\widetilde u(\mu)\|_{\Sigma} : \mu\in\mathcal{S}\}$. Here, for each new realization, we redraw the $K$ vectors $Z_1,\hdots,Z_K$, the training set $\mathcal{P}^\text{train}$, then run \cref{algo:qoi_greedy} to construct the dual reduced space $\widetilde{\mathcal{Y}}$, and finally redraw the online set $\mathcal{S}$.
We observe that for a larger tolerance $\tol$ the histograms are shifted to the left, which seems to be a bit stronger for larger $K$ (corresponding to smaller $w$).
This is due to the fact that \cref{algo:qoi_greedy} is stopped earlier and the dimension of $\widetilde{\mathcal{Y}}$ is not sufficiently large to approximate well the error estimator $\Delta(\mu)$. 
Nevertheless, we observe that the effectivity indices $\widetilde\Delta(\mu)/\|u(\mu)-\widetilde u(\mu)\|_{\Sigma}$ are always in the interval $[(\alpha w)^{-1} , (\alpha w)]$, where $\alpha \approx \tol$, as expected thanks to \cref{prop:dualErrorMultiplicative}.
This shows that, even with a rather crude approximation of the dual solutions, it is safe to use the fast-to-evaluate error estimator $\widetilde\Delta(\mu)$, as the grey area is taking into account the approximation error in the error estimator.

To guarantee that the effectivity indices lie in a user-defined interval of the form of $[c^{-1},c]$, it is sufficient to choose $\alpha$ and $w$ such that $\alpha w=c$, see \cref{prop:dualErrorMultiplicative}. As a consequence there is a degree of freedom in the choice of $\alpha$ and $w$, meaning in the choice of $\tol\approx\alpha$ and $K=K(w)$ via relation \cref{coro:truth est Condition S}. To avoid a too large shift of the histogram to the left as for instance observed for $w=2.1$ and $\tol=3$ it seems advisable to choose $\alpha$ at least as small as $w$. Additionally, the plots corresponding to $w =3.2,2.1$ and $\alpha=2$ highlight the importance of choosing $\alpha$ small enough compared to $w$ if one is interested in rather tight estimates. However, decreasing $\alpha$ has, as anticipated, a much stronger effect on the dimension of the dual reduced space (see \cref{fig:varying_alpha_w}). Therefore, it seems that for the considered test case choosing $\alpha/w \in (1/3,1)$ seems to be a good compromise between computational costs and effectivity of the error estimator. We also see that for instance $w=6.5$ and $\alpha=3$ or $\alpha=2$ yield already very good results in this direction.

\section{Conclusions}\label{sec:Conclusion}
In this paper we introduced a randomized a posteriori error estimator for low-rank approximations, which is constant-free and is both reliable and efficient at given high probability. Here, the upper and lower bound of the effectivity is chosen by the user. To derive the error estimator we exploit the concentration phenomenon of Gaussian maps. 
Exploiting the error residual relationship and approximating the associated random dual problems via projection-based model order reduction yields a fast-to-evaluate a posteriori error estimator. 
We highlight that we had to put some effort in proving the concentration inequalities but regarding the parametrized problem we only relied on its well-posedness and the definition of the adjoint operator. Therefore, there is some chance that the presented framework might be extended quite easily to \blue{more complex} problems.

To construct the dual reduced space we employed a greedy algorithm guided by a quantity of interest that assesses the quality of the fast-to-evaluate error estimator. The numerical experiments for a multi-parametric Helmholtz problem show that we obtain much smaller dual reduced spaces than with a standard greedy driven by the dual norm of the residual or with the POD. Moreover, the numerical experiments demonstrate that for moderate upper bounds for the effectivities of about $20$ the dimension of the dual reduced space needs only to be a bit more than half of the dimension of the primal reduced space. If a very tight effectivity bound of about $2$ or $3$ is desired the dual reduced spaces have to be about twice as large as the primal approximation spaces. We emphasize however that even for larger bounds of the effectivity thanks to the concentration of measure  the effectivity is still very often close to one. Furthermore, we observed only a very moderate dependence of the dimension of the dual reduced space on the number of random vectors $K$, which controls the variance of the estimator and a very mild dependence on the dimension on the (primal) reduced space. This might indicate that the error estimator will also perform well for challenging problems. Finally, we showed that to compute the fast-to-evaluate a posteriori error estimator we need to solve one dense linear system of equations of the size of the dimension of the dual reduced space. 

Due to the above the proposed a posteriori error estimator features a very favorable computational complexity and its computational costs are often about the same as the costs for the low-rank approximation or even smaller for moderate effectivity bounds. The presented error estimator can thus be more advantageous from a computational viewpoint than error estimators based on the dual norm of the residual and a (costly to estimate) stability constant or hierarchical type error estimators.

\appendix
\section{Proofs}\label{sec:proofs}
\subsection{Proof of \cref{prop:Chi2Tail}}\label{proof:Chi2Tail}

First we give a bound for $\mathbb{P}\{ Q\leq Kw^{-2} \}$. This quantity corresponds to the cumulative distribution function of the $\chi^2(K)$ distribution evaluated at $Kw^{-2}$. We have $\mathbb{P}\{ Q\leq Kw^{-2} \} = \frac{1}{\Gamma(K/2)}\gamma( \frac{K}{2} , \frac{K}{2w^2})$, where $\Gamma(\cdot)$ is the gamma function such that $\Gamma(a)=\int_0^\infty t^{a-1} e^{-t}\mathrm{d}t$ and $\gamma(\cdot,\cdot)$ the lower incomplete gamma function defined by $\gamma(a,x)=\int_0^x t^{a-1} e^{-t}\mathrm{d}t$. Following the lines of \cite{jameson_2016}, we can write $\gamma(a,x)\leq \int_0^x t^{a-1} \mathrm{d}t = \frac{1}{a}x^a$ and 
$$
 \Gamma(a) = \int_0^a t^{a-1} e^{-t}\mathrm{d}t + \int_a^\infty t^{a-1} e^{-t}\mathrm{d}t 
 \geq e^{-a}\int_0^a t^{a-1} \mathrm{d}t + a^{a-1}\int_a^\infty  e^{-t}\mathrm{d}t = 2 a^{a-1} e^{-a} ,
$$
whenever $a\geq1$. Then, if $K\geq2$ we have
\begin{equation}\label{eq:tmp6392}
 \mathbb{P}\{ Q\leq Kw^{-2} \} 
 = \frac{1}{\Gamma(K/2)} \, \gamma\Big( \frac{K}{2} , \frac{K}{2w^2}\Big) 
 \leq \frac{(K/2)e^{K/2}}{2(K/2)^{K/2}} \cdot \frac{(K/(2w^2))^{K/2}}{K/2}
 = \frac{1}{2} \Big(\frac{\sqrt{e}}{w}\Big)^{K} .
\end{equation}
Now we give a bound for $\mathbb{P}\{  Q \geq Kw^2 \}$.
Using a Markov inequality, for any $0\leq t < 1/2$ we can write
\begin{align*}
 \mathbb{P}\{  Q \geq Kw^2 \}
 &= \mathbb{P}\{ e^{ t Q } \geq e^{t Kw^2} \}
 \leq \frac{\mathbb{E}( e^{ t Q }  ) }{e^{ t Kw^2 } }
 = \frac{ (1-2t)^{-K/2} }{e^{t Kw^2 }} ,
\end{align*}
where for the last equality we used the expression for the moment-generating function of $\chi^2(K)$. The minimum of the above quantity is attained for $t=(w^2-1)/(2w^2)$ so we can write
\begin{align*}
 \mathbb{P}\{  Q \geq Kw^2 \} 
 &\leq (w^2 e^{ 1-w^2 } )^{K/2} 
 = \Big( \frac{\sqrt{e}}{w} \Big)^{K}(w^2 e^{ -w^2/2 } )^{K} \leq \Big( \frac{\sqrt{e}}{w} \Big)^{K} \frac{2^{K}}{e^{K} }
 \leq \frac{1}{2}\Big( \frac{\sqrt{e}}{w} \Big)^{K}, 
\end{align*}
for any $K\geq 3$. Together with \cref{eq:tmp6392}, the previous inequalities allows writing
$$
 \mathbb{P}\big\{  \overline{ Kw^{-2} \leq  Q \leq Kw^2 }  \big\} 
 = \mathbb{P}\{ Q\leq Kw^{-2} \} + \mathbb{P}\{  Q \geq Kw^2 \} \leq \Big( \frac{\sqrt{e}}{w} \Big)^{K} ,
$$
for any $K\geq3$, which concludes the proof.

\subsection{Proof of \cref{prop:estimate many vectors}}\label{proof:estimate many vectors}
A union bound allows writing
\begin{align*}
 &\mathbb{P}\Big\{ w^{-1}\|v\|_\Sigma \leq \| \Phi v \|_2  \leq w  \|v\|_\Sigma ~,~\forall v\in\mathcal{M} \Big\} \\
 &\quad \geq 1-\sum_{v\in\mathcal{M}} \mathbb{P}\Big\{ \overline{w^{-1}\|v\|_\Sigma \leq \| \Phi v \|_2 \leq w \|v\|_\Sigma } \Big\} \\
 &\quad= 1-(\# \mathcal{M}) ~ \mathbb{P}\big\{  \overline{ Kw^{-2} \leq  Q \leq Kw^2 } \} \geq 1-(\# \mathcal{M})  \Big( \frac{\sqrt{e}}{w} \Big)^{K},
\end{align*}
where, for the last inequality, we used Proposition \ref{prop:Chi2Tail} (assuming $w > \sqrt{e}$ and $K\geq3$ hold).
Given $0<\delta<1$, condition 
$
  K\geq \frac{\log( \# \mathcal{M}) + \log(\delta^{-1})}{\log(w/\sqrt{e})} ,
$
is equivalent to $1-(\# \mathcal{M})( \frac{\sqrt{e}}{w} )^{K} \geq 1-\delta$ and ensures that \eqref{eq:controlError} holds for all $v\in\mathcal{M}$ with probability larger than $1-\delta$. 

\subsection{Proof of \cref{prop:dualError}}\label{proof:dualError}

Let $\Psi(\mu)=K^{-1/2}[Y_1(\mu),\hdots,Y_K(\mu)]^T $ and $\widetilde\Psi(\mu)=K^{-1/2}[\widetilde Y_1(\mu),\hdots,\widetilde Y_K(\mu)]^T $ so that, from Equations \cref{eq:DualTrick} and \cref{eq: def a post est online}, we can write $\Delta(\mu) = \|\Psi(\mu) r(\mu)\|_2$ and $\widetilde\Delta(\mu) = \|\widetilde\Psi(\mu) r(\mu)\|_2$.  Using a triangle inequality we can write
$$
  | \Delta(\mu) - \widetilde\Delta(\mu) | 
  = \big| \| \Psi(\mu) r(\mu) \|_2 - \| \widetilde \Psi(\mu) r(\mu) \|_2 \big| 
  \leq \| \Psi(\mu) r(\mu) - \widetilde \Psi(\mu) r(\mu) \|_2 
$$
Dividing by $\|u(\mu)-\widetilde u(\mu)\|_\Sigma$ we can write
\begin{align*}
 \frac{| \Delta(\mu) - \widetilde\Delta(\mu) | }{\|u(\mu)-\widetilde u(\mu)\|_\Sigma}
 &\leq \frac{\| ( \Psi(\mu) - \widetilde \Psi(\mu) ) r(\mu) \|_2 }{\|u(\mu)-\widetilde u(\mu)\|_\Sigma} = \frac{\| ( \widetilde \Psi(\mu)- \Psi(\mu)) A(\mu) ( u(\mu)-\widetilde u(\mu) ) \|_2 }{\|u(\mu)-\widetilde u(\mu)\|_\Sigma} \\
 &\leq \sup_{v\in\mathbb{R}^N \backslash \{0\}} \frac{\| ( \widetilde \Psi(\mu)- \Psi(\mu)) A(\mu) v \|_2 }{\|v\|_\Sigma}  \\
 &= \sup_{\|v\|_\Sigma=1} \sqrt{ \frac{1}{K} \sum_{i=1}^K \big( ( A(\mu)^T \widetilde Y_i(\mu) - Z_i )^T v \big)^2 }\\
 &\leq \sup_{ \|v\|_\Sigma=1 } \max_{1\leq i \leq K} | ( A(\mu)^T \widetilde Y_i(\mu) - Z_i )^T v |  \\
 &= \max_{1\leq i \leq K}  \| ( A(\mu)^T \widetilde Y_i(\mu) - Z_i )^T v \|_{\Sigma^{-1}} ,
\end{align*}
which yields \cref{eq:dualError} and concludes the proof.

\subsection{Proof of \cref{cor:dualErrorAdditive}}\label{proof:dualErrorAdditive}

By \cref{prop:dualError} we have $ | \Delta(\mu) - \widetilde\Delta(\mu) | \leq \varepsilon \|u(\mu)-\widetilde u(\mu)\|_\Sigma$, which is equivalent to
$$
 \Delta(\mu) - \varepsilon \|u(\mu)-\widetilde u(\mu)\|_\Sigma \leq \widetilde\Delta(\mu) \leq \Delta(\mu) + \varepsilon \|u(\mu)-\widetilde u(\mu)\|_\Sigma.
$$
By \cref{coro:truth est S}, it holds with probability larger than $1-\delta$ that $w^{-1} \Delta(\mu) \leq \|u(\mu)-\widetilde u(\mu) \|_\Sigma \leq w \Delta(\mu)$ for all $\mu\in\mathcal{S}$. Then with the same probability we have
$$
 (w^{-1} - \varepsilon ) \|u(\mu)-\widetilde u(\mu)\|_\Sigma \leq \widetilde\Delta(\mu) \leq (w + \varepsilon )\|u(\mu)-\widetilde u(\mu)\|_\Sigma,
$$
for all $\mu\in\mathcal{S}$, which yields \cref{eq:dualErrorAdditive} and concludes the proof.

\subsection{Proof of \cref{prop:dualErrorMultiplicative}}\label{proof:dualErrorMultiplicative}
By \cref{coro:truth est S}, it holds with probability larger than $1-\delta$ that $w^{-1} \Delta(\mu) \leq \|u(\mu)-\widetilde u(\mu) \|_\Sigma \leq w \Delta(\mu)$ for all $\mu\in\mathcal{S}$. Then with the same probability we have
$$
 \|u(\mu)-\widetilde u(\mu) \| 
 \leq  w \Delta(\mu) 
 \leq  w \left( \sup_{\mu'\in\mathcal{S}} \frac{\Delta(\mu')}{\widetilde \Delta(\mu')} \right) \widetilde \Delta(\mu) 
 \overset{\cref{eq:alpha}}{\leq} (\alpha w) \widetilde \Delta(\mu),
$$
and
$$
 \|u(\mu)-\widetilde u(\mu) \| 
 \geq w^{-1} \Delta(\mu) 
 \geq w^{-1} \left( \inf_{\mu'\in\mathcal{S}} \frac{\Delta(\mu')}{\widetilde \Delta(\mu')} \right) \widetilde \Delta(\mu) 
 \overset{\cref{eq:alpha}}{\geq} (\alpha w)^{-1} \widetilde \Delta(\mu),
$$
for any $\mu\in\mathcal{S}$, which yields \cref{eq:dualErrorMultiplicative} and concludes the proof.

\subsection{Proof of \cref{prop:monolithicAlternative}}\label{proof:monolithicAlternative}

By construction, both $\widetilde Y_i(\mu)$ and $\widetilde e(\mu)$ belong to $\widetilde{\mathcal{Y}}$. Then for all $i=1,\hdots,K$ we can write
\begin{align*}
 \widetilde Y_i(\mu)^T r(\mu)
 &=\langle r(\mu), \widetilde Y_i(\mu) \rangle 
 \overset{\cref{eq:monolithicAlternative}}{=} \langle A(\mu)\widetilde e(\mu),\widetilde Y_i(\mu) \rangle \\
 & = \langle \widetilde e(\mu), A(\mu)^T \widetilde Y_i(\mu) \rangle 
 \overset{\cref{eq:monolithic}}{=} \langle \widetilde e(\mu), Z_i \rangle  
 = Z_i^T \widetilde e(\mu).
\end{align*}
Then, by definition \cref{eq: def a post est online} we can write
$$
 \widetilde \Delta(\mu) 
 = \left( \frac{1}{K} \sum_{k=1}^{K} \big(\widetilde Y_{i}(\mu)^{T} r(\mu)  \big)^{2} \right)^{1/2} 
 = \left( \frac{1}{K} \sum_{k=1}^{K} \big(Z_i^T \widetilde e(\mu)  \big)^{2} \right)^{1/2}
$$
which gives the result.

\bibliographystyle{siam}
\bibliography{random}

\end{document}